\documentclass{article}
\usepackage{arxiv}
\usepackage{amsmath,amssymb,amsfonts}
\usepackage{float}
\usepackage{algorithm,algorithmic}
\usepackage{amsthm,latexsym,float}
\usepackage[labelformat=simple]{subfig}

\usepackage{bm,multirow}        
\usepackage{graphicx}
\usepackage{xcolor}
\usepackage{placeins}
\usepackage{doi}
\graphicspath{{./figs/}}

\usepackage{hyperref}

\title{Nonintrusive model order reduction for cross-diffusion systems}

\author{ \href{https://orcid.org/0000-0003-1037-5431}{\includegraphics[scale=0.06]{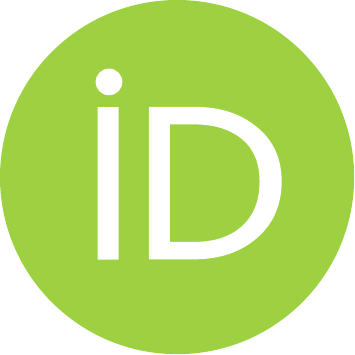}\hspace{1mm}B\"ulent Karas\"ozen} \\
    Institute of Applied Mathematics \& Department of Mathematics\\
Middle East Technical University,
 Ankara-Turkey\\
     \texttt{bulent@metu.edu.tr}\\
\And
\href{https://orcid.org/0000-0001-5262-063X}{\includegraphics[scale=0.06]{orcid.pdf}\hspace{1mm}Murat Uzunca} \\
   Department of Mathematics\\
	Sinop University,
     Sinop-Turkey \\
     \texttt{muzunca@sinop.edu.tr}\\
\And
\href{https://orcid.org/0000-0001-8952-7658}{\includegraphics[scale=0.06]{orcid.pdf}\hspace{1mm}G{\"u}den M{\"u}layim} \\
  Department of Mathematics\\
	Ad{\i}yaman University,
	Ad{\i}yaman-Turkey \\
     Institute of Applied Mathematics\\
     Middle East Technical University,  Ankara-Turkey \\
     	\texttt{gulden.mulayim@metu.edu.tr}
}

\date{}

\begin{document}
\maketitle

\begin{abstract}
In this  paper, we investigate tensor based nonintrusive reduced-order models (ROMs) for parametric cross-diffusion equations. The full-order model (FOM) consists of ordinary differential equations (ODEs) in matrix or tensor form resulting from finite-difference discretization of the differential operators by taking the advantage of Kronecker structure. The matrix/tensor differential equations are integrated in time with the implicit-explicit (IMEX) Euler method. The reduced bases, relying on a finite sample set of parameter values, are constructed in form of a two-level approach by applying higher-order singular value decomposition (HOSVD) to the space-time snapshots in tensor form, which leads to a large amount of computational and memory savings. The nonintrusive reduced approximations for an arbitrary parameter value are obtained through tensor product of the reduced basis by the parameter dependent core tensor that contains the reduced coefficients. The reduced coefficients for new parameter values are computed using radial basis function (RBF) interpolation. The efficiency of the proposed method is illustrated through numerical experiments for two-dimensional Schnakenberg and three-dimensional Brusselator cross-diffusion equations. The spatiotemporal patterns are accurately predicted by the reduced-order models with speed-up factors of orders two and three over the full-order models.
\end{abstract}


\section{Introduction}

Reaction-diffusion systems have been largely employed in literature to predict spatiotemporal patterns occurring in biological sciences, chemistry and physics. The correlation between diffusion and cross-diffusion terms may cause unstable steady solutions in form of patterns like labyrinths, spots, stripes, etc. These patterns may exhibit dynamical behavior like oscillation, annihilation, aggregation, segregation, and replication in a long time.
The common aspect of pattern formation is the interplay
between diffusion and reaction, known as diffusion-driven or Turing instability. A generalization of diffusion-driven instability is the cross-diffusion, which is characterized by a gradient in the density of one species inducing a flux in direction of another species.
Reaction-diffusion models which take into account the effects of self-diffusion as well as cross-diffusion are widely used to describe spatiotemporal dynamics of many two interacting species \cite{Amivata20,An16,Dehghan16,Gambino16,Madzvamuse15,Lin14,Zhang15,Sun12,Tulumello14}. In contrast to the classical reaction-diffusion systems without cross-diffusion, it is no longer necessary to enforce that one of the species diffuses much faster than the other for the occurrence of spatiotemporal patterns \cite{Madzvamuse15}.

Cross-diffusion systems are coupled systems of semi-linear partial differential equations (PDEs). They have been discretized in space by various methods like finite differences, finite volumes, and finite elements. In order to resolve the patterns accurately, very fine meshes in space and time are needed in numerical simulations. The effects of cross-diffusion on pattern formation in reaction-diffusion systems have been studied theoretically and numerically in many papers. Numerical simulations require fine spatial grids and long-term integration. Furthermore, multi-query simulations are required for the prediction of patterns in the parameter space.
A Cross-diffusion system involves many parameters limiting the use of
standard vector-based ordinary differential equation (ODE) solvers in time because of excessive computational costs in two and three-dimensional domains.
Under certain assumptions on the domain, one can take advantage of the Kronecker structure arising in standard space discretizations of the differential operators, and the resulting system of ODEs can be treated directly in matrix or tensor form \cite{Simoncini16,Simoncini16a,Simoncini20,Simoncini20a}.
By exploiting the structure of the diffusion matrix, the matrix/tensor based  versions of classical time integrators, such as implicit-explicit (IMEX) methods \cite{Simoncini20} allow for much finer problem discretizations.
They are based on the explicit factorization of small matrices, requiring a sequence of small matrix/tensor problems, i.e., Sylvester equations. Exploiting the spectral structure of these matrices, the computational cost is reduced further. Due to the modest size of these matrices, the computational cost per iteration can be made lower than that of the corresponding vector approaches, by working in the reduced spectral space. In this paper, we employ the matrix/tensor oriented strategy in \cite{Simoncini20} for space-time discretization of cross-diffusion systems in two and three-dimensional domains.

Simulation of the cross-diffusion systems to predict the spatiotemporal patterns for different parameter combinations take a long time and are computationally very expensive.
Reduced-order model (ROM) methods have been developed to reduce the dimension of large dynamic systems. The main idea of ROM is to construct basis functions on low-dimensional reduced space and then project onto a full-order model (FOM) to obtain reduced-order solutions.
ROMs for time-dependent parametrized PDEs have to approximate solutions as a function of time, spatial coordinates, and a parameter vector, which turns out to be more challenging.
Reduced-order modeling techniques are generally implemented in an offline-online paradigm. In the offline stage, a set of reduced
basis functions are extracted from the snapshots, i.e., a collection of high-fidelity solutions, and the reduced basis is computed by combining them. In the online phase, the FOM  is projected onto the reduced space that represents the main dynamics of the FOM, and the solutions for new parameters are computed in an efficient manner. Based on the offline-online methodology, ROM methods are classified into two categories: intrusive
and nonintrusive ROM methods. The intrusive ROM methods determine the reduced solutions by solving a reduced order model, i.e., a projection of the FOM onto the reduced space. The proper orthogonal decomposition (POD) with the Galerkin projection \cite{Berkooz93,Sirovich87} is one of the most popular tool used in intrusive ROM methods. The POD extracts the reduced basis through the singular value decomposition (SVD) of the snapshot matrix obtained by sampling in parameter space.
Then, ROM is constructed by applying Galerkin projection.
To handle this problem, some nonlinearity treatment methods are introduced, such as discrete empirical interpolation method (DEIM) \cite{Chaturantabut10}. Although all ROM methods are accurate to approximate solutions, they depend on the governing equations and discretized forms of them, that is these methods are intrusive. Another class of ROM methods is the data-driven or nonintrusive ROM (NIROM) methods which are based only on accessing to snapshots and do not use governing equations.

Contrary to a large number of papers for reduced-order modeling of patterns in fluid flows, there are few studies about the prediction of spatiotemporal patterns of reaction-diffusion equations \cite{Karasozen16,Karasozen21skt}. In this paper, we follow the NIROM approach in \cite{Audouze13,Chen18,Xiao15ns} which is based on a two-level POD approach by exploiting the matrix/tensor based  discretization of the cross-diffusion system. In the first level, the reduced basis is computed by applying the higher-order SVD (HOSVD) \cite{Lathauwer00,Vannieuwenhoven12} to the space-time snapshots related to each parameter value from a sample set of parameter values, instead of using classical SVD as in \cite{Audouze13,Chen18}. In the second level, the global set of reduced bases and coefficients of the reduced solutions are computed. The undetermined coefficients in the approximation are estimated using a nonintrusive approach based on radial basis function (RBF) approximation (in contrast to Galerkin projection).
The reduced solution for a new parameter value is then obtained by interpolating the reduced solutions with the RBF. Recently, HOSVD is used as intrusive ROM with POD \cite{Moayyedi18} and as space-time nonintrusive ROM \cite{Carlberg19}. The matrix based discretization in \cite{Simoncini20} is exploited in construction of intrusive ROMs \cite{Kirsten20} with the POD and DEIM. In \cite{Kirsten20} this approach is generalized to higher-order tensor differential equations in the framework of POD and DEIM with Galerkin projection using HOSVD. 
In this paper, the ROMs are constructed nonintrusively from space-time full-order solutions in the matrix and tensor forms with the HOSVD. The reduced solutions for new parameter values are computed by utilizing the radial basis function (RBF) interpolation.
The numerical experiments on the two-dimensional Schnakenberg and three-dimensional Brusselator cross-diffusion equations show that the patterns are predicted accurately for new parameter values. Using the nonintrusive approach with HOSVD, large amount of computer memory and computational time is saved, which is observed in high compression rates and speed-up factors of reduced-order solutions over the full-order solutions.

The rest of the paper is organized as follows. In Section~\ref{sec:fom}, we briefly describe the cross-diffusion systems and give the matrix/tensor based  space discretization by finite-differences, with the IMEX Euler time integration. The tensor based  space-time nonintrusive ROM is presented in Section~\ref{sec:rom}. Numerical results illustrating the accuracy and efficiency of the ROM methodology for the prediction of spatiotemporal patterns are given in Section~\ref{sec:num} for two examples of cross-diffusion systems: two-dimensional Schnakenberg and three-dimensional Brusselator equations. The paper ends with some conclusions in Section~\ref{sec:conc}.

\section{Full order model}
\label{sec:fom}

In this section, we briefly introduce the cross-diffusion system and describe the matrix/tensor based  discretization in space and time. We use the following notation. Scalars will be denoted as lower-case letters, vectors as bold lower-case letters, matrices are represented by capital letters or bold capital letters, and order-$d$ tensors ($d>2$) as calligraphic capital letters.

\subsection{Cross-diffusion systems}

Cross-diffusion systems are characterized by a gradient in the concentration of one species inducing a flux of another chemical species. In nature, cross-diffusion expresses the population fluxes of one species, preys,  due to the presence of the other species, predators.
The two-component cross-diffusion system is given as
\begin{equation} \label{cd}
\begin{aligned}
&u_t = d_u \nabla^2 u + d_{vu} \nabla^2 v + f(u,v),  & ({\bm x},t) \in \Omega\times (0,t_f],\\
&v_t =  d_v \nabla^2 v +  d_{uv} \nabla^2 u  + g(u,v), & ({\bm x},t) \in \Omega\times (0,t_f],\\
&\frac{\partial u}{\partial \bf{n}} = \frac{\partial v}{\partial \bf{n}}= 0,   & ({\bm x},t) \in \partial \Omega\times [0,t_f],\\
&u({\bm x},0) = u_0({\bm x}), \; v({\bm x},0)= v_0({\bm x}) & {\bm x} \in \Omega,
\end{aligned}
\end{equation}
where $\Omega \in \mathbb{R}^2$ ($\Omega \in \mathbb{R}^3$) is the spatial domain with the boundary $\partial \Omega$, ${\bm x}=(x,y)^T\in \Omega$ (${\bm x}=(x,y,z)^T\in \Omega$) is the spatial point, $\bf{n}$ is the exterior unit normal vector to the boundary, and $[0,t_f]$ is the time domain for a final time $t_f>0$.
The non-negative bounded functions $u_0({\bm x})>0$ and $v_0({\bm x})>0$ are prescribed as initial conditions. In the  cross-diffusion system \eqref{cd}, the unknown components $u({\bm x},t)$ and $v({\bm x},t)$ represent chemical concentrations or population densities.
The cross-diffusion system \eqref{cd} is a semi-linear PDE consisting of the linear diffusion parts with the Laplace operator  $\nabla^2=\partial^2/\partial x^2+\partial^2/\partial y^2$ ($\nabla^2=\partial^2/\partial x^2+\partial^2/\partial y^2+\partial^2/\partial z^2$), and nonlinear reaction terms $f(u,v)$ and $g(u,v)$.
The self-diffusion coefficients  $d_u > 0$ and $d_v>0$  are always positive whereas the cross-diffusion ones $d_{vu}$ and $d_{uv}$ can be either positive or negative. The cross-diffusion coefficient $d_{vu}$ indicates the influence of the density of $v({\bm x},t)$ to the density of $u({\bm x},t)$ so that $u({\bm x},t)$ is repelled from $v({\bm x},t)$ when $d_{vu}>0$, or otherwise $u({\bm x},t)$ is attracted to $v({\bm x},t)$ when $d_{vu}<0$. The other cross-diffusion coefficient $d_{uv}$ has the same meaning with the role of $u({\bm x},t)$ and $v({\bm x},t)$ are switched. In other words, species with positive cross-diffusion move towards the other species with the lower concentration, while in case of the negative cross-diffusion coefficients the respective species moves towards the higher concentration regime of the other species.

In the cross-diffusion systems of type \eqref{cd}, there exist a variety of reactions terms $f(u,v)$ and $g(u,v)$ with polynomial nonlinearities, such as the Schnakenberg model \cite{Gambino16,Madzvamuse15}, Brusselator model \cite{An16,Dehghan16,Lin14}, Gray-Scott model \cite{Amivata20}. Many cross-diffusion systems have nonlinear reaction terms in form of the rational functions (see for example \cite{Zhang15,Sun12,Tulumello14}).
The reaction terms in the cross-diffusion Schnakenberg model \cite{Madzvamuse15} are given by
\begin{equation} \label{Schnakenberg}
f (u, v) = \gamma(\alpha - u + u^2 v), \quad  g(u, v) = \gamma(\beta-u^2v)),
\end{equation}
where  $\gamma$ is a positive constant describing the relative strength of reaction terms.
The reaction terms of the Brusselator cross-diffusion system \cite{An16,Dehghan16,Lin14} have similar form as the one given in \eqref{Schnakenberg} of the Schnakenberg model
\begin{equation} \label{Brusselator}
f (u, v) = -(\beta + 1)u + u^2v + \alpha, \quad  g(u, v) = \beta u -u^2v.
\end{equation}

Cross-diffusion systems are also parameter dependent PDEs. In addition to the self-diffusion and cross-diffusion parameters $d_u,d_v,d_{uv},d_{vu}$, the nonlinear reaction terms includes parameters such as $\alpha,\beta$ as given in \eqref{Schnakenberg} and \eqref{Brusselator}. In this paper, we study the parameter dependent reduced-order solutions $u({\bm x},t;{\bm \mu})$ and $v({\bm x},t;{\bm \mu})$ of the system \eqref{cd}  in a parameter space ${\mathcal D}$. In this section, we suppress the parameter dependency of the states $u$ and $v$ to simplify the notation.

\subsection{Matrix and tensor based  discretization}

Semi-discretization of the cross-diffusion system \eqref{cd}  in space  with finite differences, finite elements, and spectral methods inside a hypercube in $\mathbb{R}^d$ ($d=2,3$) leads to a system of ODEs in the following form
\begin{equation} \label{cdsemv}
\begin{aligned}
\dot{\bm{u}} &= d_u A \bm{u} + d_{vu} A \bm{v}+ \bm{f}(\bm{u},\bm{v}), \quad \bm{u}(0) =\bm{u}_0, \\
\dot{\bm{v}} &= d_{uv} A \bm{u} + d_v A \bm{v} +    \bm{g}({\bm{u}},\bm{v} ), \quad \bm{v}(0) = \bm{v}_0, \\
\end{aligned}
\end{equation}
where the entries of the matrix $A$  accounts for the  spatial discretization  of  the diffusion terms including the Laplace operator $\nabla^2$ on a discrete mesh of the domain $\Omega$.    The time dependent vectors $\bm{u}(t),\bm{v}(t):[0,t_f]\rightarrow \mathbb{R}^N$ are the semi-discrete approximations of the unknown solutions $u({\bm x},t)$ and  $v({\bm x},t)$ of the system \eqref{cd}, and  $\bm{f}(\bm{u},\bm{v}),\bm{g}(\bm{u},\bm{v}):\mathbb{R}^N\times \mathbb{R}^N\rightarrow \mathbb{R}^N$ are the nonlinear vectors corresponding to the nonlinear functions $f(u,v)$ and $g(u,v)$, after spatial discretization.
All the  state vectors and nonlinear terms are evaluated componentwise at the spatial grid points. Moreover, the number $N$ denotes the degree of freedom of the discrete spatial grid. When the  finite-differences are used for the space discretization, for instance, we have $N= n_1n_2$ ($N= n_1n_2n_3$), where $n_1$ and $n_2$ ($n_1$, $n_2$ and $n_3$) are the number of the spatial nodes in $x$ and $y$-directions ($x$, $y$ and $z$-directions), respectively.

Most of the time integrators are developed for solving the semi-discretized ODEs in vector form like \eqref{cdsemv}. For an accurate simulation of the patterns of cross-diffusion systems \eqref{cd}, fine spatial discretization is required. This limits the use of standard vector-based ODE solvers in time because of the excessive computational cost and computer memory.
By exploiting the structure of the matrix of Laplace operator after space discretization, the semi-discrete ODE system  \eqref{cdsemv} can be written as a matrix/tensor differential equation. The space discretization by means of matrix/tensor based leads to the solution of linear equations with small matrices, which allows to much finer discretization of the problem and reduces the cost of full-order solutions.
In this paper, we apply  the matrix/tensor based  approach in \cite{Simoncini20} for the solution of the linear cross-diffusion systems \eqref{cd} in two and three space dimensions.

For finite differences methods, for certain finite elements techniques and spectral methods, the Laplace operator $\nabla^2$ can be discretized by means of a tensor basis. To do this, let the matrix $T_n\in\mathbb{R}^{n\times n}$ given by
$$
 T_n= \frac{1}{h^2}
\begin{pmatrix}
 -2& 2&  & & 0 \\
1& -2& 1 & &   \\
  & \ddots  & \ddots  &\ddots  &   \\
  &  &1& -2 &1 \\
 0&  &  & 2& -2
\end{pmatrix},
$$
denotes the  matrix corresponding to the discretization of the Laplace operator by centered finite differences under homogeneous Neumann boundary condition on a one-dimensional spatial mesh (an interval) $\Omega=[0,\ell]\subset\mathbb{R}$ consisting of $n$ grid points including the end points, and with the uniform mesh size $h=\ell/(n-1)$.
Then, the discretization of the Laplace operator on a rectangular domain $\Omega  = [0, \ell_x]\times [0, \ell_y]\subset\mathbb{R}^2$ leads to a matrix $A$ of the form
\begin{equation} \label{matA2}
A = I_{n_2} \otimes T_{n_1} + T_{n_2} \otimes I_{n_1} \in \mathbb{R}^{N \times N}, \qquad  N=n_1 n_2,
\end{equation}
whereas it has the form
\begin{equation} \label{matA3}
A = I_{n_3} \otimes I_{n_2} \otimes T_{n_1} + I_{n_3} \otimes  T_{n_2} \otimes I_{n_1} +  T_{n_3} \otimes I_{n_2} \otimes I_{n_1}  \in
\mathbb{R}^{N \times N \times N}, \qquad N=n_1 n_2 n_3,
\end{equation}
on a rectangular prism  $\Omega  = [0, \ell_x]\times [0, \ell_y]\times [0, \ell_z]\subset\mathbb{R}^3$. Here,  $I_{n_1}$, $I_{n_2}$ and $I_{n_3}$  are $n_1$, $n_2$ and $n_3$-dimensional identity matrices, respectively, and $\otimes$ denotes the Kronecker product. The numbers $n_1$, $n_2$ and $n_3$ denote the number of the nodes in  $x$, $y$ and $z$-directions with the mesh sizes $h_x = \ell_x/(n_1-1)$, $h_y= \ell_y/(n_2-1)$ and $h_z= \ell_z/(n_3-1)$, respectively.
Throughout the paper, we simplify the notation by taking $T_1=T_{n_1}$, $T_2=T_{n_2}$ and $T_3=T_{n_3}$ together with  $I_1=I_{n_1}$, $I_2=I_{n_2}$ and $I_3=I_{n_3}$ with the appropriate dimension.

\subsection{Full discretization on two-dimensional domains}

We consider a discrete mesh on a rectangular domain $\Omega  = [0, \ell_x]\times [0, \ell_y]\subset\mathbb{R}^2$ with the mesh sizes $h_x=\ell_x/(n_1-1)$ and $h_y=\ell_y/(n_2-1)$, and with the grid nodes ${\bm x}_{ij}=(x_i,y_j)$, where $x_i=(i-1)h_x$ and $y_j=(j-1)h_y$, $i=1,\ldots , n_1$, $j=1,\ldots , n_2$.
In order to represent the time dependent semi-discrete matrix solutions at the grid nodes, we introduce the matrix functions $U(t),V(t): [0,t_f]\mapsto  \mathbb{R}^{n_1\times n_2}$, which contain the same components of the  solution vectors $\bm{u}(t)$ and $\bm{v}(t)$ in the form $U_{ij} (t) =  u({\bm x}_{ij}, t)$ and $V_{ij} (t) =  v({\bm x}_{ij}, t)$, respectively. The rows and columns of $U$ and $V$ reflect the space discretization of the given problem in $x$ and $y$-directions, respectively.

On the other hand, the solution vectors $\bm{u}(t)$ and $\bm{v}(t)$ in the ODE system \eqref{cdsemv} can be related to the vectorization of the solution matrices $U(t)$ and $V(t)$ by the $\text{vec}(\cdot)$ operator defined by $\bm{u}(t)=\text{vec}(U(t))$ and $\bm{v}(t)=\text{vec}(V(t))$, respectively. With this operation, for instance, each column of the matrix $U(t)$ is stuck one after the other in order to obtain the vector $\text{vec}(U(t))$. This implementation satisfies a lexicographic order of the nodes in the rectangular grid for a finite difference discretization. With this notation and using the properties of the Kronecker product together with the identity \eqref{matA2}, we have that  $A{\bm{u}} = \text{vec}(T_1U +  UT_2^T)$ and $A{\bm{v}} = \text{vec}(T_1V +  VT_2^T)$. Then, the vectorial ODE system \eqref{cdsemv} can be equivalently written as the following matrix differential equation \cite{Simoncini20}
\begin{equation} \label{semimatdiff}
\begin{aligned}
\dot {U} &= d_u (T_1U + U T_2^T) + d_{vu}  (T_1 V + V T_2^T) + F(U,V), \\
\dot {V} &= d_{uv}  (T_1 U + U T_2^T) + d_v (T_1 V + V T_2^T) +  G(U,V), \\
\end{aligned}
\end{equation}
where the nonlinear matrix functions $F,G: \mathbb{R}^{n_1\times n_2}\times \mathbb{R}^{n_1\times n_2}\mapsto \mathbb{R}^{n_1\times n_2}$ are given by $F_{ij}(U,V)=f(U_{ij},V_{ij})$ and $G_{ij}(U,V)=g(U_{ij},V_{ij})$, $i=1,\ldots ,n_1$, $j=1,\ldots ,n_2$, with the property that $\bm{f}(\bm{u},\bm{v})=\text{vec}(F(U,V))$ and $\bm{g}(\bm{u},\bm{v})=\text{vec}(G(U,V))$.

Semi-discrete diffusion problems like \eqref{cdsemv} are stiff problems, which makes explicit methods inappropriate. In the presence of nonlinear reaction terms, fully implicit schemes
require a nonlinear solver, e.g., Newton, at each time step. IMEX schemes are splitting methods for ODE systems, where the stiff linear diffusion part is integrated implicitly, while the nonlinear reaction part is integrated explicitly,
as a consequence, only one linear system must be solved at each time step.
We consider the discrete times $t_k = k\Delta t$, $k = 0, \ldots ,n_t$, with the time step $ \Delta t =t_f/n_t$.  The semi-discrete matrix  differential equation \eqref{semimatdiff} is solved with IMEX Euler method \cite{Simoncini20}, which leads to the following full discrete system
\begin{equation} \label{matfom}
\begin{aligned}
\frac{U^{k+1}-U^k}{\Delta t} &= d_u (T_1 U^{k+1} + U^{k+1} T_2^T) + d_{vu}  (T_1 V^{k+1} + V^{k+1} T_2^T) + F(U^{k},V^{k}),\\
\frac{V^{k+1}-V^k}{\Delta t} &= d_{uv}  (T_1 U^{k+1} + U^{k+1} T_2^T) + d_v (T_1 V^{k+1} + V^{k+1} T_2^T) +  G(U^{k},V^{k}),\\
\end{aligned}
\end{equation}
where the full discrete solution matrices are given as $U^{k}=U(t_k)\in\mathbb{R}^{n_1\times n_2}$ and $V^{k}=V(t_k)\in\mathbb{R}^{n_1\times n_2}$, $k=0,\ldots ,n_t$, with the given initial solution matrices $U^0$ and $V^0$ satisfying $U^0_{ij} =  u_0({\bm x}_{ij})$ and $V^0_{ij} =  v_0({\bm x}_{ij})$, respectively.

The matrix/tensor formulation \eqref{semimatdiff}  has the same convergence and stability properties of the underlying time discretization methods for the classical vector differential equation \eqref{cdsemv} \cite{Simoncini20}. Exploiting the  structure of the Laplace operator in the linear part and using finer grids, highly accurate full-order solutions are obtained and  the computational cost is much reduced.
After collecting the alike terms in the full discrete system \eqref{matfom}, we obtain a system of linear matrix equation in the form of the Sylvester equation
\begin{equation} \label{sylvester}
\begin{aligned}
(I_1-d_u \Delta t T_1)U^{k+1} - d_u \Delta t U^{k+1}T_2^T -  d_{vu} \Delta t T_1 V^{k+1} - d_{vu}\Delta t V^{k+1}T_2^T &= U^k + \Delta tF(U^k,V^k),\\
(I_1-d_v \Delta t T_1)V^{k+1} - d_v \Delta t V^{k+1}T_2^T -  d_{uv} \Delta t T_1 U^{k+1} - d_{uv}\Delta t U^{k+1}T_2^T &= V^k + \Delta tG(U^k,V^k),
\end{aligned}
\end{equation}
which is solved for the solution matrices $U^{k+1}$ and $V^{k+1}$.

The matrix/tensor methods can be made more efficient by computing a-priori spectral decomposition of the coefficient matrices/tensors of not too large sizes \cite{Simoncini20}. Assuming that the matrices $T_1$ and $T_2^T$ are diagonalizable, the solution of the Sylvester equation \eqref{sylvester} is accelerated by the use of the eigenvalue decomposition of the matrices $T_1$ and $T_2^T$.
Let the eigenvalue decompositions $T_1=X\Lambda^{(x)} X^{-1}$ and $T_2^T=Y\Lambda^{(y)} Y^{-1}$ are given, with the matrices $X\in\mathbb{R}^{n_1\times n_1}$ and $Y\in\mathbb{R}^{n_2\times n_2}$ of nonsingular vectors, and the diagonal matrices $\Lambda^{(1)} = \text{diag}(\lambda_1^{(1)},\ldots,\lambda_{n_1}^{(1)})$ and $\Lambda^{(2)} = \text{diag}(\lambda_1^{(2)},\ldots,\lambda_{n_2}^{(2)})$ of the eigenvalues $\lambda_i^{(1)}$ and $\lambda_j^{(2)}$, $i=1,\ldots, n_1$, $j=1,\ldots, n_2$.
Multiplying both equations in \eqref{sylvester} from left by $X^{-1}$ and from right by $Y$, substituting eigenvalue decompositions $T_1=X\Lambda^{(1)} X^{-1}$ and $T_2^T=Y\Lambda^{(2)} Y^{-1}$, and setting $\hat{U}^k = X^{-1} U^{k} Y$ and $\hat{V}^k = X^{-1} V^{k} Y$, we reach the system of matrix equations
\begin{equation}\label{sylvester2}
\begin{aligned}
(I_1-d_u \Delta t \Lambda^{(1)})\hat{U}^{k+1} - d_u \Delta t \hat{U}^{k+1} \Lambda^{(2)} - d_{vu} \Delta t \Lambda^{(1)} \hat{V}^{k+1}  - d_{vu} \Delta t \hat{V}^{k+1}\Lambda^{(2)} &= Q_1^k, \\
(I_1-d_v \Delta t \Lambda^{(1)})\hat{V}^{k+1} - d_v \Delta t \hat{V}^{k+1} \Lambda^{(2)} - d_{uv} \Delta t \Lambda^{(1)} \hat{U}^{k+1}  - d_{uv} \Delta t \hat{U}^{k+1}\Lambda^{(2)} &= Q_2^k, \\
\end{aligned}
\end{equation}
where
\begin{equation}\label{q12}
\begin{aligned}
Q_1^k &=  X^{-1} (U^k + \Delta tF(U^k,V^k)) Y,\\
Q_2^k &=  X^{-1} (V^k + \Delta tG(U^k,V^k)) Y.
\end{aligned}
\end{equation}

The matrix equation \eqref{sylvester2} in which all the coefficient matrices on the left hand sides are diagonal matrices, can be easily solved componentwise. Therewith, the entries of the solution matrices $U^{k+1},V^{k+1} \in \mathbb{R}^{n_1 \times n_2}$ are given  as the following $2\times2$  linear system of equations
\begin{equation} \label{matrixform}
\begin{aligned}
\begin{pmatrix}
S^{11}_{ij} & S^{12}_{ij} \\
S^{21}_{ij} & S^{22}_{ij} \\
\end{pmatrix}
\begin{pmatrix}
\hat{U}^{k+1}_{ij} \\ \hat{V}^{k+1}_{ij}
\end{pmatrix} =
\begin{pmatrix}
(Q_1^k)_{ij} \\ (Q_2^k)_{ij}\\
\end{pmatrix}, \quad i=1,\ldots n_1, \; j=1,\ldots n_2,
\end{aligned}
\end{equation}
where for  $p,q=1,2$, the entries of the matrices $S^{pq} \in \mathbb{R}^{n_1 \times n_2}$ are given by
\begin{equation*}
\begin{aligned}
S^{11}_{ij} &= 1-d_u \Delta t (\lambda_i^{(1)} + \lambda_j^{(2)}), \qquad S^{12}_{ij} &= - d_{vu} \Delta t (\lambda_i^{(1)} + \lambda_j^{(2)}), \\
S^{22}_{ij} &= 1-d_v \Delta t (\lambda_i^{(1)} + \lambda_j^{(2)}), \qquad S^{21}_{ij} &= - d_{uv} \Delta t (\lambda_i^{(1)} + \lambda_j^{(2)}). \\
\end{aligned}
\end{equation*}
The solution of the $2\times2$  linear system \eqref{matrixform} for fixed $i$ and $j$ can  be written as
\begin{equation*}
\begin{aligned}
\begin{pmatrix}
\hat{U}^k_{ij} \\ \hat{V}^n_{ij}
\end{pmatrix} =
\begin{pmatrix}
S^{11}_{ij} & S^{12}_{ij} \\
S^{21}_{ij} & S^{22}_{ij}
\end{pmatrix}^{-1}
\begin{pmatrix}
(Q_1^k)_{ij} \\ (Q_2^k)_{ij}
\end{pmatrix},
\end{aligned}
\end{equation*}
where for the determinant $\mid S_{ij}\mid=S^{11}_{ij} S^{22}_{ij} - S^{12}_{ij}S^{21}_{ij}$, the $2 \times 2$ inverse matrix can be  easily calculated as
$$
\begin{pmatrix}
	S^{11}_{ij} & S^{12}_{ij} \\
	S^{21}_{ij} & S^{22}_{ij} \\
\end{pmatrix}^{-1}
=
\frac{1}{ \mid S_{ij}\mid }
\begin{pmatrix}
S^{22}_{ij} & - S^{12}_{ij} \\
-S^{21}_{ij} & S^{11}_{ij} \\
\end{pmatrix},
$$
Finally, by introducing the matrices $L^{pq} \in \mathbb{R}^{n_1 \times n_2}$, $p,q=1,2$, with the entries
\begin{equation} \label{l12}
L^{11}_{ij}= \frac{S^{22}_{ij}}{\mid S_{ij}\mid}, \quad L^{12}_{ij}= \frac{-S^{12}_{ij}}{\mid S_{ij}\mid}, \quad L^{21}_{ij}= \frac{-S^{21}_{ij}}{\mid S_{ij}\mid}, \quad L^{22}_{ij}= \frac{S^{11}_{ij}}{\mid S_{ij}\mid},
\end{equation}
the solution of the $2\times2$  linear system \eqref{matrixform} are given by
$$
\hat{U}^k_{ij} = L^{11}_{ij} (Q_1^k)_{ij} + L^{12}_{ij} (Q_2^k)_{ij}, \quad \hat{V}^k_{ij} = L^{21}_{ij} (Q_1^k)_{ij} + L^{22}_{ij} (Q_2^k)_{ij}.
$$
The unknown solution matrices $U^{k+1}$ and $V^{k+1}$ can then be recovered by projecting back as $U^{k+1}=X\hat{U}^kY^{-1}$ and $V^{k+1}=X\hat{V}^kY^{-1}$, which can be written in terms of the matrices $L^{pq}$, $Q_1^k$ and $Q_2^k$ as
\begin{align*}
U^{k+1} &= X(L^{11} \odot  Q_1^k + L^{12} \odot Q_2^k)Y^{-1},\\
V^{k+1} &= X(L^{21} \odot Q_1^k + L^{22} \odot  Q_2^k)Y^{-1},
\end{align*}
where $\odot$ denotes the Hadamard (element by element) product.
The solution process to compute the full discrete solution matrices $U^{k+1}$ and $V^{k+1}$ by using IMEX Euler method applied to the semi-discrete linear matrix differential equation \eqref{semimatdiff}, and by utilizing the eigenvalue decompositions $T_1=X\Lambda^{(1)} X^{-1}$ and $T_2^T=Y\Lambda^{(2)} Y^{-1}$ is given in Algorithm~\ref{alg}.

\begin{algorithm}
\caption{Solution process on a single time step \label{alg}}
\begin{algorithmic}[1]
	\STATE\textbf{Input:} Known solution matrices $U^k$ and $V^k$, eigenvectors $X$ and $Y$, eigenvalues $\Lambda^{(1)}$ and $\Lambda^{(2)}$
	\STATE\textbf{Output:} Unknown solution matrices $U^{k+1}$ and $V^{k+1}$
	\STATE \quad Compute the matrices $Q_1^k$ and $Q_2^k$ from \eqref{q12}
	\STATE \quad Compute the matrices $L^{11}$, $L^{12}$, $L^{21}$ and $L^{22}$ from \eqref{l12}
	\STATE \quad Compute the solution matrices $U^{k+1}$ and $V^{k+1}$ as
	\begin{align*}
U^{k+1} &= X(L^{11} \odot  Q_1^k + L^{12} \odot Q_2^k)Y^{-1}\\
V^{k+1} &= X(L^{21} \odot Q_1^k + L^{22} \odot  Q_2^k )Y^{-1}
\end{align*}
\end{algorithmic}
\end{algorithm}

The overall computational cost of solving two and three dimensional cross-diffusion systems is drastically reduced using the matrix/tensor formulation with the explicit-implicit time integration and using spectral decomposition. The computation of the spectral decomposition of the matrices $T_1$ and $T_2$ are performed once at the beginning of the integration.

\subsection{Full discretization on three-dimensional domains}

We consider a discrete mesh on a rectangular prism $\Omega  = [0, \ell_x]\times [0, \ell_y]\times [0, \ell_z]\subset\mathbb{R}^3$ with the mesh sizes $h_x=\ell_x/(n_1-1)$, $h_y=\ell_y/(n_2-1)$ and $h_z=\ell_z/(n_3-1)$, and with the grid nodes ${\bm x}_{ijl}=(x_i,y_j,z_l)$, where $x_i=(i-1)h_x$, $y_j=(j-1)h_y$ and $z_l=(l-1)h_z$, $i=1,\ldots , n_1$, $j=1,\ldots , n_2$, $l=1,\ldots , n_3$.

The matrix oriented approach in \cite{Simoncini20} can be extended to the cross-diffusion systems \eqref{cd} for the three-dimensional case following \cite{Simoncini16}.
At each time step $t \in [0, t_f ]$,  let $U_{i,j}^l (t) =  u({\bm x}_{ijl}, t)$ and $V_{i,j}^l (t) =  v({\bm x}_{ijl}, t)$  denote  the approximate semi-discrete solutions at the grid nodes ${\bm x}_{ijl}$. Then, we introduce tall matrix functions $\bm{U}(t), \bm{V}(t):[0,t_f]\mapsto \mathbb{R}^{(n_1n_3) \times n_2}$ defined as
$$
\bm{U}(t) = \left[
\begin{matrix}
U^1 (t)  \\
U^2 (t) \\
\vdots \\
U^{n_3} (t)
\end{matrix}\right],\qquad
\bm{V}(t) = \left[
\begin{matrix}
V^1 (t)  \\
V^2 (t) \\
\vdots \\
V^{n_3} (t)
\end{matrix}\right].
$$
With this notation and an appropriate ordering of the nodes by the $\text{vec}(\cdot)$ operation introduced before, the terms in the system \eqref{cdsemv} including the matrix $A\in\mathbb{R}^{(n_1n_2n_3)\times (n_1n_2n_3)}$ and the solution vectors $\bm{u},\bm{v}: [0,t_f]\mapsto \mathbb{R}^{(n_1n_2n_3)}$ can be  written as
\begin{equation} \label{semimatdiff3}
\begin{aligned}
A\bm{u} &= \text{vec}\left( (I_3\otimes T_1) \bm{U} + \bm{U} T_2^T  + (T_3\otimes I_1) \bm{U} \right)  \\
A\bm{v} &= \text{vec}\left( (I_3\otimes T_1) \bm{V} + \bm{V} T_2^T  + (T_3\otimes I_1) \bm{V}\right) \\
\end{aligned}
\end{equation}
At the discrete times $t_k = k\Delta t$, $k = 0, \ldots ,n_t$, let $\bm{U}^k=\bm{U}(t_k)\in\mathbb{R}^{(n_1n_3)\times n_2}$ and $\bm{V}^k=\bm{V}(t_k)\in\mathbb{R}^{(n_1n_3)\times n_2}$  denote the full discrete solution matrices at the time $t_k$. Then, using the identity \eqref{semimatdiff3}, application of the IMEX Euler method, similar to the two-dimensional case, yields the following Sylvester equation as a matrix differential equation
\begin{equation} \label{sylvester3}
\begin{aligned}
(I_{13}-d_u \Delta t \widehat{T})\bm{U}^{k+1} - d_u \Delta t \bm{U}^{k+1}T_2^T -  d_{vu} \Delta t \widehat{T} \bm{V}^{k+1} - d_{vu}\Delta t \bm{V}^{k+1}T_2^T &= \bm{U}^k + \Delta t\bm{F}(\bm{U}^k,\bm{V}^k),\\
(I_{13}-d_v \Delta t \widehat{T})\bm{V}^{k+1} - d_v \Delta t \bm{V}^{k+1}T_2^T -  d_{uv} \Delta t \widehat{T} \bm{U}^{k+1} - d_{uv}\Delta t \bm{U}^{k+1}T_2^T &= \bm{V}^k + \Delta t\bm{G}(\bm{U}^k,\bm{V}^k),
\end{aligned}
\end{equation}
where $\widehat{T}:=T_1\oplus T_3=(I_3\otimes T_1+T_3\otimes I_1)\in\mathbb{R}^{(n_1n_3)\times (n_1n_3)}$ with $\oplus$ denoting the Kronecker sum, and $I_{13}$ is the identity matrix of size $(n_1n_3)$.
The above Sylvester equation can be solved similar to the two-dimensional case. Here, it needs only the use of the eigenvalue decomposition of the matrix $\widehat{T}=\widehat{X}\widehat{\Lambda}\widehat{X}^{-1}$ in place of the eigenvalue decomposition of the matrix $T_1$. However, the square matrix $\widehat{T}$ is of dimension $(n_1n_3)$ which makes inefficient the computation of the eigenvalue decomposition of $\widehat{T}$. Instead, we use the eigenvalue decomposition of the matrices $T_1$, $T_2^T$ and $T_3$ of smaller size, and we use the properties of Kronecker sum. Let the eigenvalue decompositions $T_1=X\Lambda^{(1)} X^{-1}$, $T_2^T=Y\Lambda^{(2)} Y^{-1}$ and $T_3=Z\Lambda^{(3)} Z^{-1}$ are given, with the matrices $X\in\mathbb{R}^{n_1\times n_1}$, $Y\in\mathbb{R}^{n_2\times n_2}$ and $Z\in\mathbb{R}^{n_3\times n_3}$ of nonsingular vectors and the diagonal matrices $\Lambda^{(1)} = \text{diag}(\lambda_1^{(1)},\ldots,\lambda_{n_1}^{(1)})$, $\Lambda^{(2)} = \text{diag}(\lambda_1^{(2)},\ldots,\lambda_{n_2}^{(2)})$ and $\Lambda^{(3)} = \text{diag}(\lambda_1^{(3)},\ldots,\lambda_{n_3}^{(3)})$ of the eigenvalues $\lambda_i^{(1)}$, $\lambda_j^{(2)}$ and $\lambda_l^{(3)}$, $i=1,\ldots, n_1$, $j=1,\ldots, n_2$, $l=1,\ldots, n_3$.
Then, by the use of the properties of the Kronecker sum, the eigenvalue decomposition of the matrix $\widehat{T}$ with the nonsingular vector $\widehat{X}$ and the diagonal matrix of the eigenvalues $\widehat{\Lambda} = \text{diag}(\widehat{\lambda}_1,,\ldots ,\widehat{\lambda}_{n_1n_3})$ are given by
$$
\widehat{X}=(Z\otimes I_1 )(I_3\otimes X) ,  \quad \widehat{\Lambda}= (\Lambda^{(3)}\otimes I_1 ) + (I_3\otimes \Lambda^{(1)}).
$$
Similar to the two-dimensional case,  the Sylvester equation \eqref{sylvester3} can be efficiently solved through multiplying both the equations in \eqref{sylvester3} from left by $\widehat{X}^{-1}$ and from right by $Y$, substituting  $\widehat{T}=\widehat{X}\widehat{\Lambda}\widehat{X}^{-1}$ and $T_2^T=Y\Lambda^{(2)} Y^{-1}$, and setting $\hat{\bm{U}}^k = \widehat{X}^{-1} \bm{U}^{k} Y$ and $\hat{\bm{V}}^k = \widehat{X}^{-1} \bm{V}^{k} Y$, yielding the system
\begin{equation}\label{sylvester31}
\begin{aligned}
(I_{13}-d_u \Delta t \widehat{\Lambda})\hat{\bm{U}}^{k+1} - d_u \Delta t \hat{\bm{U}}^{k+1} \Lambda^{(2)} - d_{vu} \Delta t \widehat{\Lambda} \hat{\bm{V}}^{k+1}  - d_{vu} \Delta t \hat{\bm{V}}^{k+1}\Lambda^{(2)} &= \bm{Q}_1^k, \\
(I_{13}-d_v \Delta t \widehat{\Lambda})\hat{\bm{V}}^{k+1} - d_v \Delta t \hat{\bm{V}}^{k+1} \Lambda^{(2)} - d_{uv} \Delta t \widehat{\Lambda} \hat{\bm{U}}^{k+1}  - d_{uv} \Delta t \hat{\bm{U}}^{k+1}\Lambda^{(2)} &= \bm{Q}_2^k, \\
\end{aligned}
\end{equation}
where all the coefficient matrices on the left hand sides are again diagonal matrices, and the right hand side matrices are given by
\begin{equation*}
\begin{aligned}
\bm{Q}_1^k &=  \widehat{X}^{-1} (\bm{U}^k + \Delta t\bm{F}(\bm{U}^k,\bm{V}^k)) Y,\\
\bm{Q}_2^k &=  \widehat{X}^{-1} (\bm{V}^k + \Delta t\bm{G}(\bm{U}^k,\bm{V}^k)) Y.
\end{aligned}
\end{equation*}

\section{Nonintrusive reduced-order model}
\label{sec:rom}

In this section, we consider the following parametrized form of the vectorial cross diffusion system \eqref{cdsemv}
\begin{equation} \label{cdsemvpar}
\begin{aligned}
\dot{\bm{u}}^{\theta} &= d_u A \bm{u}^{\theta} + d_{vu} A \bm{v}^{\theta} + \bm{f}(\bm{u}^{\theta},\bm{v}^{\theta};\theta), \quad \bm{u}^{\theta}(0) =\bm{u}^{\theta}_0, \\
\dot{\bm{v}}^{\theta} &= d_{uv} A \bm{u}^{\theta} + d_v A \bm{v}^{\theta} +    \bm{g}({\bm{u}}^{\theta},\bm{v}^{\theta};{\theta} ), \quad \bm{v}^{\theta}(0) = \bm{v}^{\theta}_0, \\
\end{aligned}
\end{equation}
where the superscript ${\theta}\in\textsc{P}$ indicates the parameter dependency of the solutions, and $\textsc{P}$ is a  set of admissible values of the parameter $\theta$ which may stand for either parameter in the system. In most cases, the system \eqref{cdsemvpar} needs to be solved several times by the value of parameter $\theta$ differs. In this paper, by using a finite training set $\textsc{P}_T=\{\theta_1,\ldots ,\theta_{n_p}\}\subset\textsc{P}$, we aim to construct a nonintrusive ROM to the system \eqref{cdsemvpar} in order to cheaply obtain approximate solutions for a given parameter value $\theta\in\textsc{P}$ not necessarily from the training set $\textsc{P}_T$, i.e., $\theta\notin\textsc{P}_T$.

Reduced-order modelling methodology commonly relies on a data obtained from either an experiment or solutions of a discrete system like \eqref{cdsemvpar}, which is named as snapshot data. In our case, solving the parametrized cross diffusion system \eqref{cdsemvpar} through either the matrix system \eqref{sylvester} for a two-dimensional domain ($d=2$) or the matrix system \eqref{sylvester3} for a three-dimensional domain ($d=3$), and with a suitable arrangement of the dimensions, we can obtain a set of snapshots $\{{\mathcal U}^{\theta}(t_k)\}_{k=1}^{n_t}$ and $\{{\mathcal V}^{\theta}(t_k)\}_{k=1}^{n_t}$ in the form of an order-$d$ tensor (multidimensional array) with ${\mathcal U}^{\theta}(t_k)\in \mathbb{R}^{n_1\times \cdots\times n_d}$ and ${\mathcal V}^{\theta}(t_k)\in \mathbb{R}^{n_1\times \cdots\times n_d}$.
Here, each dimension of the tensors ${\mathcal U}^{\theta}(t_k)$ and ${\mathcal V}^{\theta}(t_k)$ corresponds to one of the respective spatial directions, for instance, with $d=3$, ${\mathcal U}^{\theta}_{i_1i_2i_3}(t_k)=u^{\theta}({\bm x}_{i_1i_2i_3}, t_k)$, $i_s=1,\ldots ,n_s$, $s=1,2,3$. Then, we form the following order-$(d+1)$ tensors of snapshot data related to a given parameter $\theta$
\begin{equation} \label{snap}
\begin{aligned}
&{\mathcal X}^{\theta,u}\in \mathbb{R}^{n_1\times \cdots\times n_d\times n_{d+1}} , & {\mathcal X}^{\theta,u}_{i_1\cdots i_di_{d+1}}={\mathcal U}^{\theta}_{i_1\cdots i_d}(t_{i_{d+1}}),\\
&{\mathcal X}^{\theta,v}\in \mathbb{R}^{n_1\times \cdots\times n_d\times n_{d+1}} , & {\mathcal X}^{\theta,v}_{i_1\cdots i_di_{d+1}}={\mathcal V}^{\theta}_{i_1\cdots i_d}(t_{i_{d+1}}),
\end{aligned}
\end{equation}
where for easy notation we set the size of the final dimension related to the time as $n_{d+1}:=n_t$.

The standard POD approach to construct the reduced basis solutions for many training parameter values is costly.
The two-level POD, known also as nested POD, is often used in ROM applications for parametrized PDEs to reduce the computational cost of constructing the spatial and temporal basis functions \cite{Audouze13,Chen18}. Usually, a snapshot data which is in the form of columns consisting of the solutions in vector form, is used in this ROM methodology.
In the first level, a set of POD basis are computed for the snapshot data related to each parameter $\theta_i\in\textsc{P}_T$. Then, in the second level, a global POD basis is constructed by applying SVD to the set of POD basis computed in the first level \cite{Audouze13,Chen18}. Finally, the space-time coefficients are determined using RBF in a nonintrusive way without resorting to Galerkin projection.

For the snapshot data in the form of an order-$2$ tensor, i.e., a matrix, the POD basis in the first level are computed by applying SVD or eigenvalue decomposition to the snapshot matrix. However, in our case, each snapshot data given in \eqref{snap} is an order-$(d+1)$ tensor with $d=2,3$.
Here, in the first level of nested POD, we compute the POD modes using HOSVD of the snapshot tensors ${\mathcal X}^{\theta,u}$ and ${\mathcal X}^{\theta,v}$. The HOSVD is a favorite algorithm for computing low-rank approximation of the Tucker decomposition of a tensor \cite{Tucker66}. In the following, we will describe HOSVD to compute the POD modes of the snapshot tensor ${\mathcal X}^{\theta,u}$, the POD modes of the snapshot tensor ${\mathcal X}^{\theta,v}$ can be computed similarly.
Like any multidimensional array, the order-$(d+1)$ snapshot tensor ${\mathcal X}^{\theta,u}\in \mathbb{R}^{n_1\times \cdots\times n_d\times n_{d+1}}$ admits the following Tucker decomposition \cite{Vannieuwenhoven12,Tucker66,Minster20}
$$
{\mathcal X}^{\theta,u}_{{i_1\cdots i_di_{d+1}}} = \sum_{j_1=1}^{n_1}\cdots \sum_{j_d=1}^{n_d}\sum_{j_{d+1}=1}^{n_{d+1}} {\mathcal S}^{\theta,u}_{{j_1\cdots j_dj_{d+1}}}\Phi^{(1),\theta,u}_{{i_1j_1}}\cdots\Phi^{(d),\theta,u}_{{i_dj_d}}\Phi^{(d+1),\theta,u}_{{i_{d+1}j_{d+1}}},
$$
or in a suitable compact form
\begin{equation}\label{tucker}
{\mathcal X}^{\theta,u} = \left( \Phi^{(1),\theta,u},\ldots ,\Phi^{(d),\theta,u},\Phi^{(d+1),\theta,u} \right) \cdot {\mathcal S}^{\theta,u},
\end{equation}
where the order-$(d+1)$ tensor ${\mathcal S}^{\theta,u}\in \mathbb{R}^{n_1\times \cdots\times n_d\times n_{d+1}}$ is called the core tensor, and each orthonormal matrix $\Phi^{(j),\theta,u}\in \mathbb{R}^{n_j\times n_j}$,  $j=1,\ldots ,(d+1)$, is called a factor matrix. In other words, a tensor can be decomposed into a core tensor that is multiplied by a matrix along each mode, which are orthonormal and can be viewed as the principal components of each modes. The HOSVD aims firstly to compute the factor matrices $\Phi^{(j),\theta,u}$, $j=1,\ldots,(d+1)$. This process is done by applying SVD to the mode-$j$ unfolding (matricization) ${\mathcal X}^{\theta,u}_{(j)}$ of the tensor ${\mathcal X}^{\theta,u}$, where a  mode-$j$ unfolding ${\mathcal X}^{\theta,u}_{(j)}$ is a matrix of size $n_j\times \prod_{i\neq j}n_i$, and its columns are mode-$j$ fibers of the tensor ${\mathcal X}^{\theta,u}$ \cite{Vannieuwenhoven12,Tucker66,Minster20}. Then, the $j$th factor matrix is given by the left singular vectors of the mode-$j$ unfolding of the tensor ${\mathcal X}^{\theta,u}$
$$
{\mathcal X}^{\theta,u}_{(j)} = \Phi^{(j),\theta,u}\Sigma^{(j),\theta,u} \left( \psi^{(j),\theta,u}\right)^T, \quad j=1,\ldots , (d+1),
$$
where the diagonal matrix $\Sigma^{(j),\theta,u}\in \mathbb{R}^{n_j\times n_j}$ includes on its diagonal elements the singular values $\sigma_i\left({\mathcal X}^{\theta,u}_{(j)}\right)\geq 0$ of the mode-$j$ unfolding ${\mathcal X}^{\theta,u}_{(j)}$, $i=1,\ldots ,n_j$. After computation of the factor matrices, the core tensor can be calculated as
$$
{\mathcal S}^{\theta,u} = \left( \left(\Phi^{(1),\theta,u}\right)^T,\ldots ,\left(\Phi^{(d),\theta,u}\right)^T,\left(\Phi^{(d+1),\theta,u}\right)^T \right) \cdot {\mathcal X}^{\theta,u}.
$$
By the use of HOSVD, the factor matrices take place of the POD modes required in the first level of the nested POD, each of which corresponds to one of the either space direction or temporal dimension. In addition, the space-time coefficients are contained in the core tensor ${\mathcal S}^{\theta,u}$, therefore there is no need to determine them through the solution of a separate system as in \cite{Audouze13,Chen18}.

On the other hand, the decomposition \eqref{tucker} does not provide a low-rank approximation yet, it requires $\prod_{i=1}^{d+1}n_i + \sum_{i=1}^{d+1}n_i^2$ numbers to be stored. The HOSVD can be employed to construct a low multilinear rank approximation to a tensor, where it provides a compressed representation in the Tucker decomposition. One approach is the truncated HOSVD (T-HOSVD) which was first introduced in \cite{Lathauwer00}. The T-HOSVD algorithm aims to compute each factor matrix separately, and it relies on truncating each mode-$j$ unfolding ${\mathcal X}^{\theta,u}_{(j)}$ of the snapshot tensor ${\mathcal X}^{\theta,u}$ according to a given truncation criteria or an a priori given target ranks for each dimension $j$. Let for some target rank $\bm{r}^{\theta ,u}=(r^{\theta ,u}_1,\ldots,r^{\theta ,u}_{d+1})$ with $r^{\theta ,u}_j<n_j$ for each $j=1,\ldots ,(d+1)$, the truncated SVD of the mode-$j$ unfolding ${\mathcal X}^{\theta,u}_{(j)}$ is given by
\begin{equation}\label{redsvd}
\begin{aligned}
{\mathcal X}^{\theta,u}_{(j)} &= \Phi^{(j),\theta,u}\Sigma^{(j),\theta,u} \left( \psi^{(j),\theta,u}\right)^T =
\begin{bmatrix}
\bar{\Phi}^{(j),\theta,u} & \tilde{\Phi}^{(j),\theta,u}
\end{bmatrix}
\begin{bmatrix}
\bar{\Sigma}^{(j),\theta,u} & \\
              & \tilde{\Sigma}^{(j),\theta,u}
\end{bmatrix}
\begin{bmatrix}
\left( \bar{\psi}^{(j),\theta,u}\right)^T\\
\left( \tilde{\psi}^{(j),\theta,u}\right)^T
\end{bmatrix},
\end{aligned}
\end{equation}
where $\bar{\Phi}^{(j),\theta,u}\in \mathbb{R}^{n_j \times r^{\theta ,u}_j}$ contains  the first $r^{\theta ,u}_j$ left singular vectors from $\Phi^{(j),\theta,u}$, retained singular values are contained in $\bar{\Sigma}^{(j),\theta,u}\in \mathbb{R}^{r^{\theta ,u}_j \times r^{\theta ,u}_j}$, and $\tilde{\Sigma}^{(j),\theta,u}$ contains the truncated singular values. Using the truncated factor matrices $\bar{\Phi}^{(j),\theta,u}$, we can calculate the truncated (reduced) core tensor $\bar{{\mathcal S}}^{\theta,u}\in \mathbb{R}^{r^{\theta ,u}_1\times \cdots\times r^{\theta ,u}_d\times r^{\theta ,u}_{d+1}}$ using the formula
$$
\bar{{\mathcal S}}^{\theta,u} = \left( \left(\bar{\Phi}^{(1),\theta,u}\right)^T,\ldots ,\left(\bar{\Phi}^{(d),\theta,u}\right)^T,\left(\bar{\Phi}^{(d+1),\theta,u}\right)^T \right) \cdot {\mathcal X}^{\theta,u}.
$$
Then, a rank-$(r^{\theta ,u}_1,\dots, r^{\theta ,u}_d,r^{\theta ,u}_{d+1})$ approximation $\bar{{\mathcal X}}^{\theta,u}\in \mathbb{R}^{n_1\times \cdots\times n_d\times n_{d+1}}$ to the snapshot tensor ${\mathcal X}^{\theta,u}\in \mathbb{R}^{n_1\times \cdots\times n_d\times n_{d+1}}$ can be obtained as
\begin{equation}\label{thosvd}
{\mathcal X}^{\theta,u} \approx \bar{{\mathcal X}}^{\theta,u} = \left( \bar{\Phi}^{(1),\theta,u},\ldots ,\bar{\Phi}^{(d),\theta,u},\bar{\Phi}^{(d+1),\theta,u} \right) \cdot \bar{{\mathcal S}}^{\theta,u},
\end{equation}
where it stores only $\prod_{i=1}^{d+1}r^{\theta ,u}_i + \sum_{i=1}^{d+1}n_ir^{\theta ,u}_i$ numbers. The quantity that to what extend the memory saving is obtained, can be visualized by the following compression factor \cite{Lorento10}
\begin{equation}\label{comp}
C_F = \frac{\prod_{i=1}^{d+1}n_i + \sum_{i=1}^{d+1}n_i^2}{\prod_{i=1}^{d+1}r^{\theta ,u}_i + \sum_{i=1}^{d+1}n_ir^{\theta ,u}_i},
\end{equation}
which gives the saved memory in percentage by the formula $100(1-1/C_F)$. The larger the compression factor $C_F$ the much more the memory is saved.

Although, T-HOSVD provides a low multilinear rank approximation, the SVD computations of the unfoldings may be expensive, since the same full rank tensor is used to obtain each unfolding.
Another approach to construct a low multilinear rank approximation is the sequentially truncated HOSVD (ST-HOSVD) \cite{Vannieuwenhoven12,Minster20}, which is a variation of the usual T-HOSVD. In the ST-HOSVD, instead of throwing away most of the work performed by each SVD, SVD is performed sequentially on a reduced tensor along
all dimensions. Starting from the initial core tensor $\bar{{\mathcal S}}^{\theta,u,(0)}:={\mathcal X}^{\theta,u}$, ST-HOSVD computes a sequence of core tensors $\bar{{\mathcal S}}^{\theta,u,(j)}$ to reach the reduced core tensor $\bar{{\mathcal S}}^{\theta,u}=\bar{{\mathcal S}}^{\theta,u,(d+1)}$ following an ordering $\bm{p}=(p_1,\ldots,p_{d+1})$ which is a permutation of the index set $(1,2,\ldots ,d+1)$. In the $j$th stage, the truncated factor matrix $\bar{\Phi}^{(p_j),\theta,u}$ of the mode-$p_j$ unfolding of the core tensor $\bar{{\mathcal S}}^{\theta,u,(j-1)}$ is computed, and the new core tensor $\bar{{\mathcal S}}^{\theta,u,(j)}$ is calculated by projecting the previous one onto the subspace spanned by the columns of the computed factor matrix $\bar{\Phi}^{(p_j),\theta,u}$. The ST-HOSVD algorithm is given in Algorithm~\ref{alg:sthosvd} \cite{Minster20}.

\begin{algorithm}
	\caption{ST-HOSVD for tensor data of $u$ component\label{alg:sthosvd}}
	\begin{algorithmic}[1]
		\STATE \textbf{Input:} Snapshot tensor ${\mathcal X}^{\theta,u}$, processing order $\bm{p}=(p_1,\ldots ,p_{d+1})$, target ranks $\bm{r}^{\theta ,u}=(r^{\theta ,u}_1,\ldots ,r^{\theta ,u}_{d+1})$
		\STATE \textbf{Output:} Truncated factor matrices $\{\bar{\Phi}^{(1),\theta,u},\ldots ,\bar{\Phi}^{(d+1),\theta,u}\}$
		\STATE Set $\bar{{\mathcal S}}^{\theta,u,(0)}={\mathcal X}^{\theta,u}$
		\FOR{$j = 1$ to $d+1$}
		\STATE Obtain mode-$p_j$ unfolding $\bar{{\mathcal S}}^{\theta,u,(j-1)}_{(p_j)}$
		\STATE Apply truncated SVD to the unfolding $\bar{{\mathcal S}}^{\theta,u,(j-1)}_{(p_j)}$ for target rank $r^{\theta ,u}_{p_j}$
		\STATE Get the factor matrix $\bar{\Phi}^{(p_j),\theta,u}$
		\STATE Update the unfolding $\bar{{\mathcal S}}^{\theta,u,(j-1)}_{(p_j)} \; \leftarrow \; \bar{\Sigma}^{(p_j),\theta,u}\left( \bar{\psi}^{(p_j),\theta,u}\right)^T$
		\STATE Obtain updated core tensor $\bar{{\mathcal S}}^{\theta,u,(j)} \; \leftarrow \; \bar{{\mathcal S}}^{\theta,u,(j-1)}_{(p_j)}$ in tensor form
		\ENDFOR
	\end{algorithmic}
\end{algorithm}

Unlike T-HOSVD, the process in the ST-HOSVD is sequential, therefore the order in which the modes are processed affects the accuracy of the approximation and the speed of the process.
In \cite{Vannieuwenhoven12}, a heuristic is proposed that attempts to minimize the number of operations required to compute the dominant subspace. Processing first the dimension with the lowest size may even reduce the rank of the remaining terms, i.e., $n_{p_1}\leq n_{p_2}\leq \cdots \leq n_{p_{d+1}}$. In this way, more energy is forced into fewer modes.
Computing the T-HOSVD can be more expensive than the ST-HOSVD, while the ST-HOSVD requires fewer floating point operations to compute the approximation. Although, T-HOSVD and ST-HOSVD approximations may differ in accuracy for an ordering $\bm{p}\neq (1,\ldots,d+1)$,
both T-HOSVD and ST-HOSVD approximations satisfy the same error bounds \cite{Vannieuwenhoven12}
\begin{equation}\label{froberr}
\min_{1\leq j\leq d+1} \|\tilde{\Sigma}^{(j),\theta,u} \|_F^2 \leq \|{\mathcal X}^{\theta,u} - \bar{{\mathcal X}}^{\theta,u}\|_F^2 \leq \sum_{j=1}^{d+1}  \|\tilde{\Sigma}^{(j),\theta,u} \|_F^2,
\end{equation}
where $\|\cdot\|_F$ is the Frobenius norm, and $\tilde{\Sigma}^{(j),\theta,u}$ contains the truncated singular values given in \eqref{redsvd}.

In order to construct the nonintrusive ROM through the nested POD, we first form the set of snapshot tensors $\{{\mathcal X}^{\theta_i,u}\}_{i=1}^{n_p}$ and $\{{\mathcal X}^{\theta_i,v}\}_{i=1}^{n_p}$ from the solutions of the parametrized cross diffusion system related to each sample parameter $\theta_i\in \textsc{P}_T$.
Then, in the first level of the nested POD,  we apply ST-HOSVD to the snapshot tensors ${\mathcal X}^{\theta_i,u}$ and ${\mathcal X}^{\theta_i,v}$, and we collect related to each sample parameter $\theta_i\in \textsc{P}_T$, the truncated factor matrices $\{\bar{\Phi}^{(j),\theta_i,u}\}_{i=1}^{n_p}$ and $\{\bar{\Phi}^{(j),\theta_i,v}\}_{i=1}^{n_p}$ with the target ranks $\bm{r}^{\theta_i,u}=(r^{\theta_i,u}_1,\ldots ,r^{\theta_i,u}_{d+1})$ and $\bm{r}^{\theta_i,v}=(r^{\theta_i,v}_1,\ldots ,r^{\theta_i,v}_{d+1})$, respectively, $j=1,\ldots , (d+1)$.
Then, in the second level of the nested POD, we compute the truncated global factor matrices $\widehat{\Phi}^{(j),u}\in \mathbb{R}^{n_j\times \widehat{r}^u_j}$ and $\widehat{\Phi}^{(j),v}\in \mathbb{R}^{n_j\times \widehat{r}^v_j}$, $j=1,\ldots ,(d+1)$, as the truncated left singular vectors obtained by the application of the truncated SVD to the collections $\bar{\Phi}^{(j),u}$ and $\bar{\Phi}^{(j),v}$ of the factor matrices defined by
\begin{align*}
\bar{\Phi}^{(j),u} &= [\bar{\Phi}^{(j),\theta_1,u} \; \bar{\Phi}^{(j),\theta_2,u} \; \cdots \; \bar{\Phi}^{(j),\theta_{n_p},u}] \in \mathbb{R}^{n_j\times \bar{r}^u_j}, \\
\bar{\Phi}^{(j),v} &=  [\bar{\Phi}^{(j),\theta_1,v} \; \bar{\Phi}^{(j),\theta_2,v} \; \cdots \; \bar{\Phi}^{(j),\theta_{n_p},v}] \in \mathbb{R}^{n_j\times \bar{r}^v_j},
\end{align*}
where the numbers $\bar{r}^u_j:=r^{\theta_1,u}_j+\cdots +r^{\theta_{n_p},u}_j$ and $\bar{r}^v_j:=r^{\theta_1,v}_j+\cdots +r^{\theta_{n_p},v}_j$ denote the column size of the collection of the factor matrices, which are the sum of the target ranks $r^{\theta_i,u}_j$ and $r^{\theta_i,v}_j$ of each unfolding ${\mathcal X}^{\theta_i,u}_{(j)}$ and ${\mathcal X}^{\theta_i,v}_{(j)}$, respectively, $i=1,\ldots ,n_p$, $j=1,\ldots ,(d+1)$.
In addition, the numbers $\widehat{r}^u_j<\bar{r}^u_j$ and $\widehat{r}^v_j<\bar{r}^v_j$ are the target ranks of the collections $\bar{\Phi}^{(j),u}$ and $\bar{\Phi}^{(j),v}$ of the factor matrices, respectively.
Note that the truncated global factor matrices $\widehat{\Phi}^{(j),u}$ and $\widehat{\Phi}^{(j),v}$ are independent of the parameter $\theta$, they rely on the sample parameter set $\textsc{P}_T$.
Once the truncated global factor matrices $\widehat{\Phi}^{(j),u}$ and $\widehat{\Phi}^{(j),v}$ are obtained, we compute for each sample parameter $\theta_i\in \textsc{P}_T$, the core tensors $\widehat{{\mathcal S}}^{\theta_i,u}\in \mathbb{R}^{\widehat{r}^u_1\times \cdots\times \widehat{r}^u_d\times \widehat{r}^u_{d+1}}$ and $\widehat{{\mathcal S}}^{\theta_i,v}\in \mathbb{R}^{\widehat{r}^v_1\times \cdots\times \widehat{r}^v_d\times \widehat{r}^v_{d+1}}$ as
\begin{align*}
\widehat{{\mathcal S}}^{\theta_i,u} &= \left( \left(\widehat{\Phi}^{(1),u}\right)^T,\ldots ,\left(\widehat{\Phi}^{(d),u}\right)^T,\left(\widehat{\Phi}^{(d+1),u}\right)^T \right) \cdot {\mathcal X}^{\theta_i,u}, & \\
\widehat{{\mathcal S}}^{\theta_i,v} &= \left( \left(\widehat{\Phi}^{(1),v}\right)^T,\ldots ,\left(\widehat{\Phi}^{(d),v}\right)^T,\left(\widehat{\Phi}^{(d+1),v}\right)^T \right) \cdot {\mathcal X}^{\theta_i,v}, & i=1,\ldots ,n_p.
\end{align*}

Finally, the parameter dependent nonintrusive ROM solution tensors $\widehat{\mathcal X}^{u}(\theta ),\widehat{\mathcal X}^{v}(\theta )\in \mathbb{R}^{n_1 \times \cdots \times n_d\times n_{d+1}}$ for an arbitrary parameter $\theta\in \textsc{P}$ can be efficiently obtained by the formulas
\begin{equation} \label{nirom}
\begin{aligned}
{\mathcal X}^{\theta,u} \approx \widehat{\mathcal X}^{u}(\theta ) &= \left( \widehat{\Phi}^{(1),u},\ldots ,\widehat{\Phi}^{(d),u},\widehat{\Phi}^{(d+1),u}\right) \cdot \widehat{\mathcal S}^{u}(\theta ),\\
{\mathcal X}^{\theta,v} \approx \widehat{\mathcal X}^{v}(\theta ) &= \left( \widehat{\Phi}^{(1),v},\ldots ,\widehat{\Phi}^{(d),v},\widehat{\Phi}^{(d+1),v}\right) \cdot \widehat{\mathcal S}^{v}(\theta ),
\end{aligned}
\end{equation}
where the parameter dependent global core tensors $\widehat{\mathcal S}^{u}(\theta ): \textsc{P}\mapsto \mathbb{R}^{\widehat{r}^u_1\times \cdots\times \widehat{r}^u_d\times \widehat{r}^u_{d+1}}$ and $\widehat{\mathcal S}^{v}(\theta ): \textsc{P}\mapsto \mathbb{R}^{\widehat{r}^v_1\times \cdots\times \widehat{r}^v_d\times \widehat{r}^v_{d+1}}$ stands for the data of undetermined coefficients, which can be easily determined by a variety of methods. Here, each entry of the global core tensors are expanded using RBFs as follows
\begin{equation} \label{rbf}
\begin{aligned}
\widehat{\mathcal S}^{u}_{j_1\cdots j_dj_{d+1}}(\theta ) &= \sum_{j=1}^{n_p} \gamma_{j_1\cdots j_dj_{d+1}}^{u,j} \varPsi(\omega_j(\theta )), \\
\widehat{\mathcal S}^{v}_{j_1\cdots j_dj_{d+1}}(\theta ) &= \sum_{j=1}^{n_p} \gamma_{j_1\cdots j_dj_{d+1}}^{v,j} \varPsi(\omega_j(\theta )),
\end{aligned}
\end{equation}
where $\varPsi(\omega)$ denote the radial basis kernel function with $\omega_j(\theta ) = \mid \theta-\theta_j\mid$, and the scalars $\gamma_{j_1\cdots j_dj_{d+1}}^{u,j}$ and $\gamma_{j_1\cdots j_dj_{d+1}}^{v,j}$ are the coefficients to be determined. RBF is a real-valued function whose value depends on the distance from center point so that $\varPsi(\omega) = \varPsi(\|\omega\|)$ is a radial function. There exist well-known RBFs. Here, Gaussian RBF $\varPsi(\omega)=e^{(-\omega^2/2\rho)}$ is used, where the parameter $\rho$, in our case, is given by $\rho = (\max_i{\theta_i}-\min_i{\theta_i})/n_p$.

In order to compute the undetermined coefficients $\gamma_{j_1\cdots j_dj_{d+1}}^{u,j}$ and $\gamma_{j_1\cdots j_dj_{d+1}}^{v,j}$, we use the core tensors $\widehat{{\mathcal S}}^{\theta_i,u}$ and $\widehat{{\mathcal S}}^{\theta_i,v}$.
Setting $\theta = \theta_{i}$ in \eqref{rbf} with the properties that $\widehat{\mathcal S}^{u}(\theta_i)=\widehat{\mathcal S}^{\theta_i,u}$ and $\widehat{\mathcal S}^{v}(\theta_i)=\widehat{\mathcal S}^{\theta_i,v}$, we obtain that
\begin{align*}
\widehat{\mathcal S}^{u}_{j_1\cdots j_dj_{d+1}}(\theta_i ) &= \widehat{\mathcal S}^{\theta_i,u}_{j_1\cdots j_dj_{d+1}} = \sum_{j=1}^{n_p} \gamma_{j_1\cdots j_dj_{d+1}}^{u,j} \varPsi(\omega_j(\theta _i)), \\
\widehat{\mathcal S}^{v}_{j_1\cdots j_dj_{d+1}}(\theta_i ) &= \widehat{\mathcal S}^{\theta_i,v}_{j_1\cdots j_dj_{d+1}} = \sum_{j=1}^{n_p} \gamma_{j_1\cdots j_dj_{d+1}}^{v,j} \varPsi(\omega_j(\theta _i)), & i= 1,...,n_p,
\end{align*}
which leads to the following linear systems of equations
\begin{equation*}
\sum_{j=1}^{n_p} B_{i j} \gamma_{j_1\cdots j_dj_{d+1}}^{u,j} = \widehat{\mathcal S}^{\theta_i,u}_{j_1\cdots j_dj_{d+1}}, \quad \sum_{j=1}^{n_p} B_{i j} \gamma_{j_1\cdots j_dj_{d+1}}^{v,j} = \widehat{\mathcal S}^{\theta_i,v}_{j_1\cdots j_dj_{d+1}},\quad i= 1,...,n_p,
\end{equation*}
where the entries of the symmetric interpolation matrix $B\in\mathbb{R}^{n_p\times n_p}$ is given by $B_{i j}=\varPsi(\omega_j(\theta _i))$. In short, for a given new parameter value $\theta \in\textsc{P}$, once the undetermined coefficients $\gamma_{j_1\cdots j_dj_{d+1}}^{u,j}$ and $\gamma_{j_1\cdots j_dj_{d+1}}^{v,j}$ in \eqref{rbf} are computed by the RBF interpolation, the nonintrusive ROM solution tensors $\widehat{\mathcal X}^{u}(\theta )$ and $\widehat{\mathcal X}^{v}(\theta )$ in \eqref{nirom} are calculated, by which the nonintrusive ROM approximations $\widehat{u}({\bm x}, t; \theta)\approx u({\bm x}, t; \theta)$ and $\widehat{v}({\bm x}, t; \theta)\approx v({\bm x}, t; \theta)$ can be cheaply obtained as
$$
\widehat{u}({\bm x}_{i_1\cdots i_d}, t_k; \theta) = \widehat{\mathcal X}_{i_1\cdots i_dk}^{u}(\theta ), \qquad \widehat{v}({\bm x}_{i_1\cdots i_d}, t_k; \theta) = \widehat{\mathcal X}_{i_1\cdots i_dk}^{v}(\theta ).
$$

\section{Numerical results}
\label{sec:num}

In this section we report about the numerical tests for the two-dimensional Schnakenberg \eqref{Schnakenberg} and three-dimensional Brusselator \eqref{Brusselator} cross-diffusion systems. All the simulations are performed on a machine with Intel CoreTM i7 2.5 GHz 64 bit CPU, 8 GB RAM, Windows 10, using 64 bit MatLab R2014.
For both problems, the initial conditions are taken as random periodic perturbation around the equilibrium solutions $u_e$ and $v_e$
\begin{itemize}
	\item Schnakenberg \cite{Madzvamuse15}:
	$$u_0({\bm x}) = u_e +\text{rand}({\bm x})/100, \quad v_0({\bm x}) = v_e +\text{rand}({\bm x})/100,$$
	with $(u_e,v_e)=(0.55,0.9917)$,
	\item Brusselator \cite{Lin14}:
	$$u_0({\bm x}) = u_e +\text{rand}({\bm x})/3, \quad v_0({\bm x}) = v_e +\text{rand}({\bm x})/10,$$
	with $(u_e,v_e)=(6,0.1667)$,
\end{itemize}
where $\text{rand}({\bm x})$ is the MatLab's random function producing a multi-dimensional array of the same dimension as ${\bm x}$, with the entries are uniformly distributed random numbers between $0$ and $1$.

For a given parameter value $\theta\in\textsc{P}$, the accuracy of the corresponding reduced approximations are measured using the time averaged relative errors
\begin{equation}\label{relerr}
\begin{aligned}
\|u-\widehat{u}\|_{\text{rel}}&=\frac{1}{n_t}\sum_{k=1}^{n_t}\frac{\|{\mathcal U}^{\theta}(t_k)-\widehat{{\mathcal U}}^{\theta}(t_k)\|_{F}}{\|{\mathcal U}^{\theta}(t_k)\|_{F}}, \\
\|v-\widehat{v}\|_{\text{rel}}&=\frac{1}{n_t}\sum_{k=1}^{n_t}\frac{\|{\mathcal V}^{\theta}(t_k)-\widehat{{\mathcal V}}^{\theta}(t_k)\|_{F}}{\|{\mathcal V}^{\theta}(t_k)\|_{F}},
\end{aligned}
\end{equation}
where ${\mathcal U}^{\theta}(t),{\mathcal V}^{\theta}(t):[0,t_f]\mapsto \mathbb{R}^{n_1\times \cdots \times n_d}$ are the discrete FOM solutions while  $\widehat{{\mathcal U}}^{\theta}(t),\widehat{{\mathcal V}}^{\theta}(t):[0,t_f]\mapsto \mathbb{R}^{n_1\times \cdots \times n_d}$ are the discrete ROM approximations in the form of order-$d$ tensor $(d=2,3)$.

In both examples, in order to obtain the truncated factor matrices in the first level of nested POD, say $\bar{\Phi}^{(j),\theta,u}$, the target rank $r^{\theta,u}_j$ for each unfolding ${\mathcal X}^{\theta,u}_{(j)}$, $j=1,\ldots ,(d+1)$, is determined compatibly with the error bound \eqref{froberr} in Frobenius norm, so that the following criteria is satisfied for a user given tolerance $\tau_1>0$
\begin{equation}\label{cr1}
\frac{\sqrt{\sum_{i=r^{\theta,u}_j+1}^{n_j}\sigma_i\left({\mathcal X}^{\theta,u}_{(j)}\right)}}{\sqrt{\sum_{i=1}^{n_j}\sigma_i\left({\mathcal X}^{\theta,u}_{(j)}\right)}}< \tau_1.
\end{equation}
On the other hand, in the second level of nested POD, in order obtain the truncated global factor matrices, say $\widehat{\Phi}^{(j),u}$, we apply the truncated SVD to the collection
$
\bar{\Phi}^{(j),u}=[\bar{\Phi}^{(j),\theta_1,u} \; \bar{\Phi}^{(j),\theta_2,u} \; \cdots \; \bar{\Phi}^{(j),\theta_{n_p},u}]\in \mathbb{R}^{n_j\times \bar{r}^u_j}
$
($\bar{r}^u_j=r^{\theta_1,u}_j+\cdots +r^{\theta_{n_p},u}_j$ ) of the truncated factor matrices
with the target rank $\widehat{r}^u_j<\bar{r}^u_j$ which is determined so that the following energy criteria is satisfied for a user given tolerance $\tau_2>0$
\begin{equation}\label{cr2}
\frac{\sum_{i=1}^{\widehat{r}^u_j}\sigma_i\left(\bar{\Phi}^{(j),u}\right)}{\sum_{i=1}^{\bar{r}^u_j}\sigma_i\left(\bar{\Phi}^{(j),u}\right)}\geq 1- \tau_2,
\end{equation}
where $\sigma_i\left(\bar{\Phi}^{(j),u}\right)\geq 0$ are the singular values of the collection $\bar{\Phi}^{(j),u}$ of the truncated factor matrices computed in the first level. In our simulations, we choose the user given tolerances scaling as $\tau_1\sim 10^{-2}$ and $\tau_2\sim 10^{-8}$.

\subsection{Schnakenberg equation}

We consider the two-dimensional Schnakenberg model \cite{Madzvamuse15}
in  the square domain $\Omega = [0,0.5]^2\subset\mathbb{R}^2$. We solve the problem through the matrix equation \eqref{sylvester2} with the mesh sizes $h_x=h_y=0.005$ for the number of grid points $n_1=n_2=101$. The final time is taken as $t_f=5$ with the time step size $\Delta t=0.001$, leading to the third dimension $n_3=5001$ of the snapshot tensors related to the time. For this problem,
we fix the system parameters  $d_u=d_v=d_{vu}=1$, $\gamma=200$, $\alpha =0.25$, $\beta =0.3$, and vary the cross-diffusion parameter $\theta :=d_{uv}$ in the set of admissible values $\textsc{P}=[0.4,0.8]$. As the finite training set of parameter $\theta$, we take the values (including the boundary values) uniformly distributed on $\textsc{P}$ with the increment $0.1$, i.e., $\textsc{P}_T=\{0.4,0.5,0.6,0.7,0.8\}$ with the number of sample parameter values $n_p=5$.

\begin{figure}[H]
	\centering
	\includegraphics[width=0.8\textwidth]{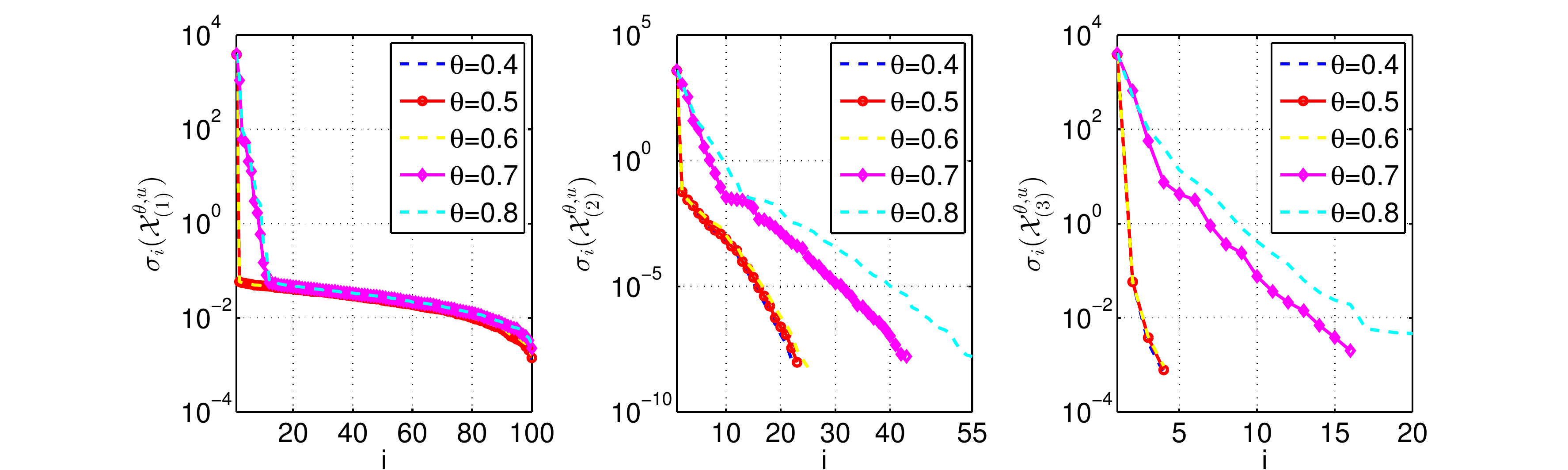}	
	\includegraphics[width=0.8\textwidth]{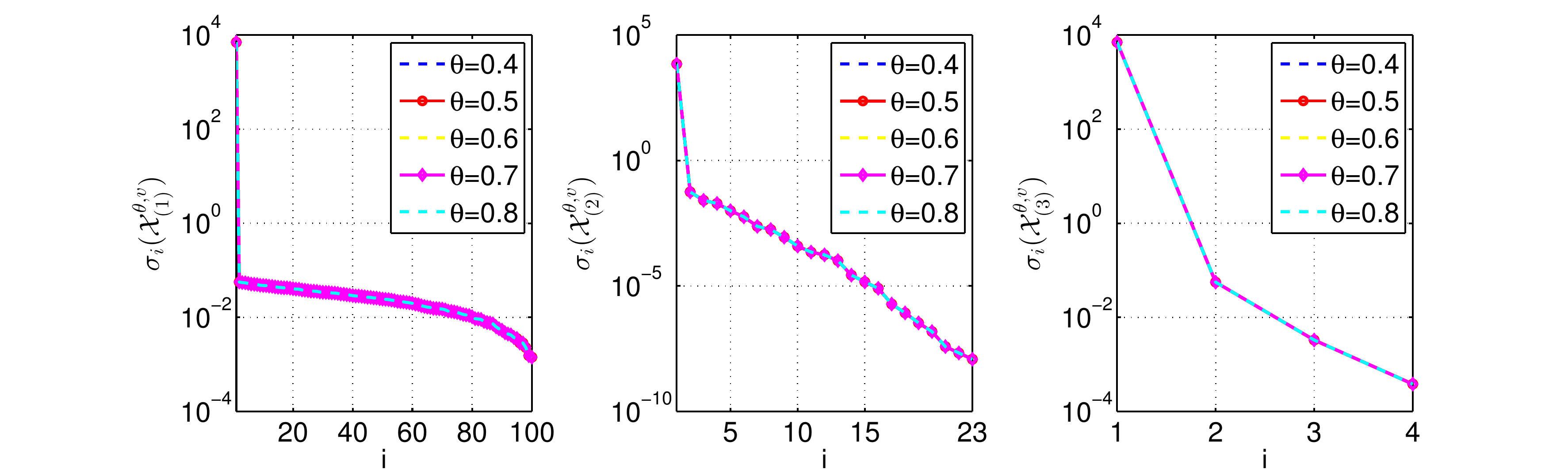}
	\caption{Schnakenberg model: Level I singular values of unfoldings ${\mathcal X}^{\theta,u}_{(j)}$ (top) and ${\mathcal X}^{\theta,v}_{(j)}$ (bottom) \label{sing2d}}
\end{figure}

\begin{figure}[htb!]
	\centering
	\includegraphics[width=0.6\textwidth]{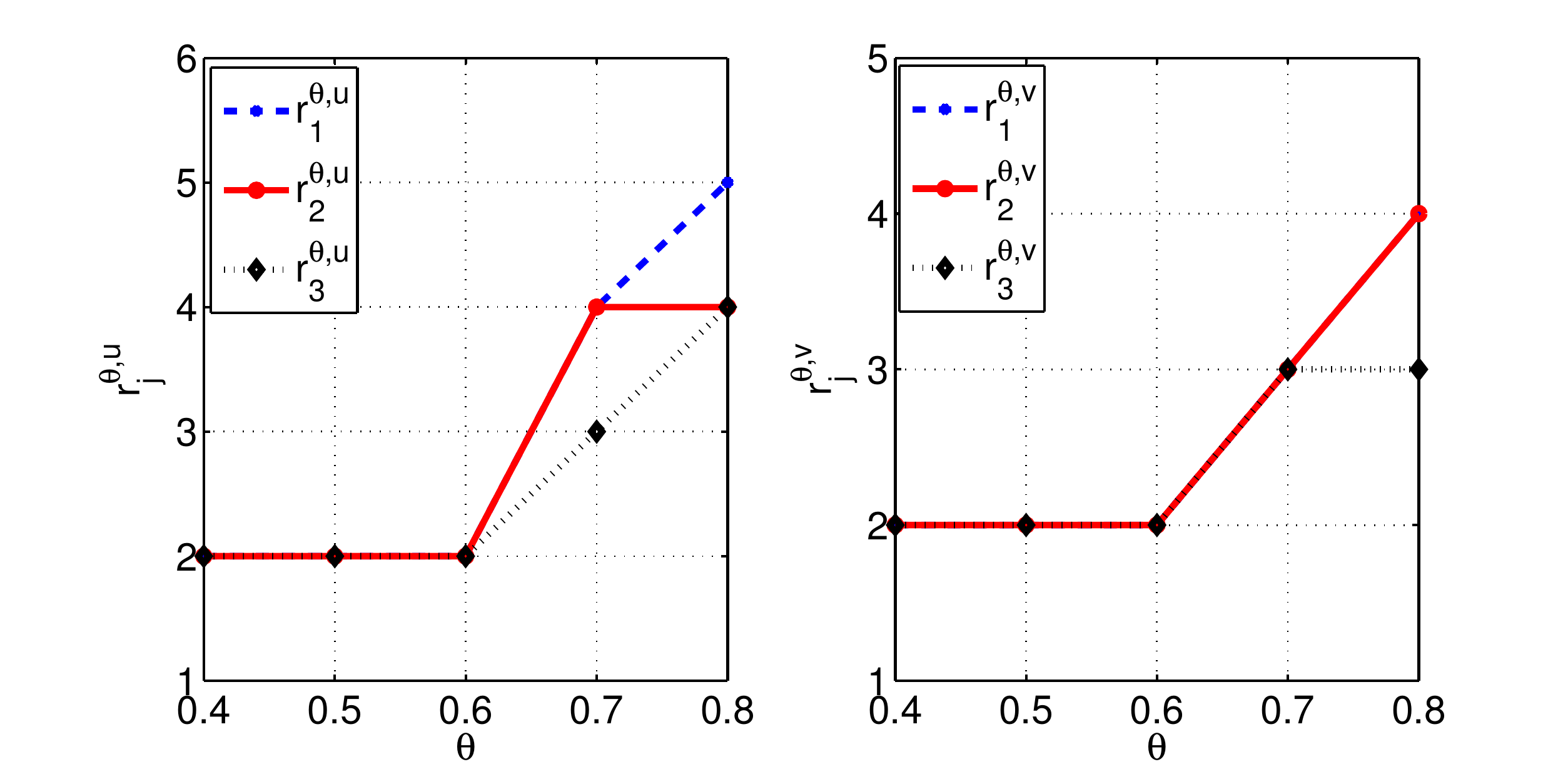}	
	\caption{Schnakenberg model: Sample parameter values vs target ranks for unfoldings ${\mathcal X}^{\theta,u}_{(j)}$ (left) and ${\mathcal X}^{\theta,v}_{(j)}$ (right) at Level I \label{numr2d}}
\end{figure}

In \figurename~\ref{sing2d}, we give the decay of the singular values $\sigma_i\left({\mathcal X}^{\theta,u}_{(j)}\right)$ and $\sigma_i\left({\mathcal X}^{\theta,v}_{(j)}\right)$ of the unfoldings ${\mathcal X}^{\theta,u}_{(j)}$ and ${\mathcal X}^{\theta,v}_{(j)}$ of the order-$3$ snapshot tensors ${\mathcal X}^{\theta,u}$ and ${\mathcal X}^{\theta,v}$, respectively, $j=1,2,3$, related to each sample parameter value $\theta\in\textsc{P}_T$. According to the criteria \eqref{cr1}, the computed target ranks $r^{\theta,u}_j$ and $r^{\theta,v}_j$ required by the HOSVD in the first level of nested POD are presented in \figurename~\ref{numr2d}, which shows in accordance with the singular values in \figurename~\ref{sing2d} that enough energetic part of the unfoldings are recovered.

\begin{figure}[htb!]
	\centering
	\includegraphics[width=0.7\textwidth]{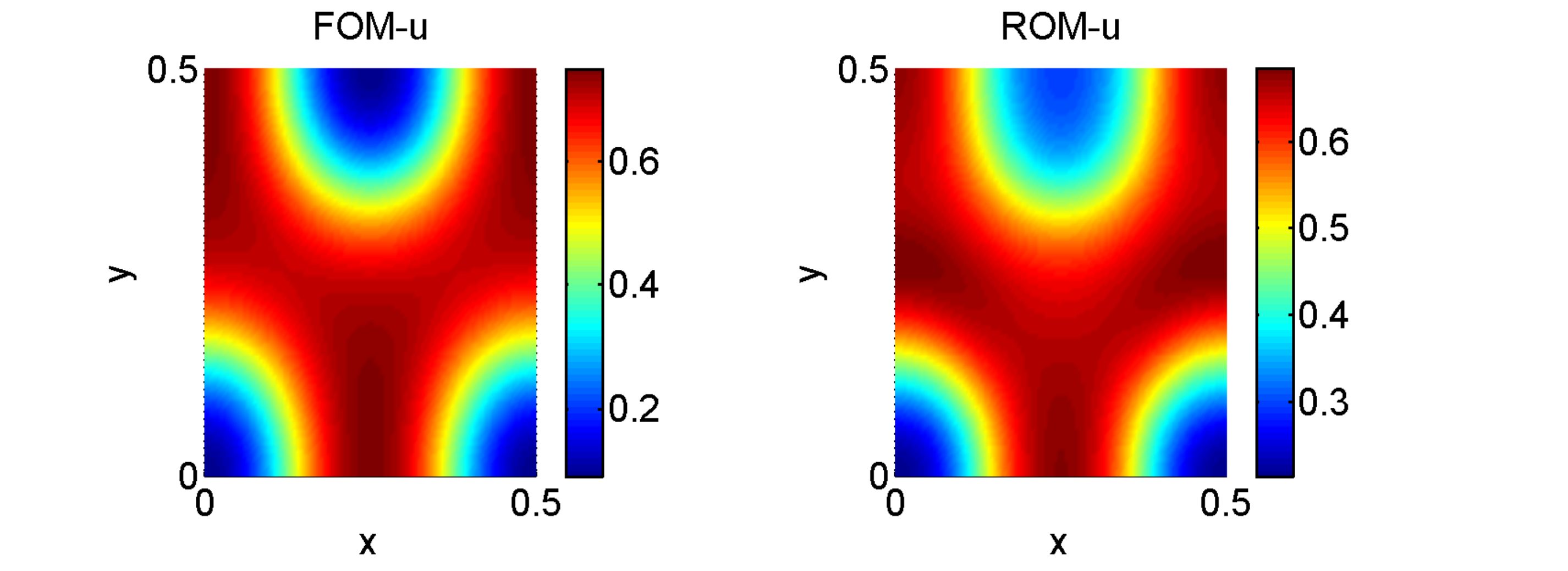}	
	\includegraphics[width=0.7\textwidth]{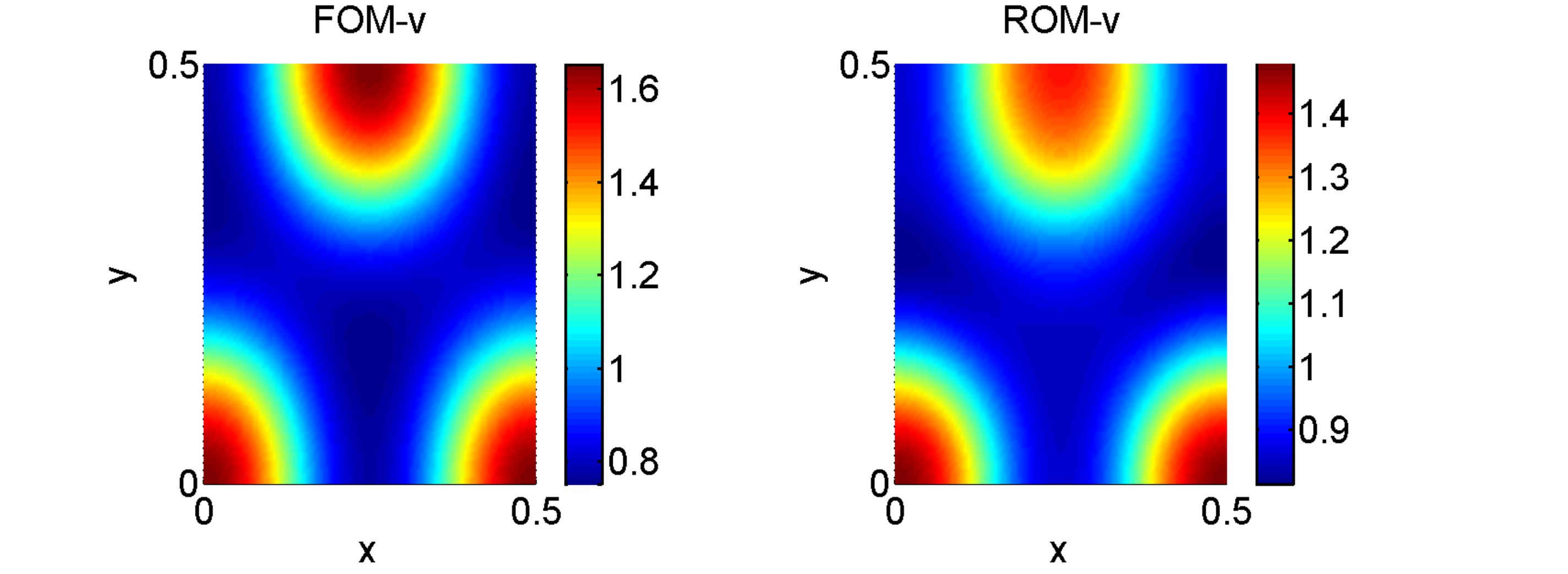}
	\caption{Schnakenberg model: FOM and ROM patterns at the final time for $\theta =0.65$\label{plot2d}}
\end{figure}

The FOM solutions $u({\bm x}, t; \theta)$ and $v({\bm x}, t; \theta)$ together with the nonintrusive ROM approximations $\widehat{u}({\bm x}, t; \theta)$ and $\widehat{v}({\bm x}, t; \theta)$ at the final time $t_f=5$ for the parameter value $\theta =0.65\notin\textsc{P}_T$ are given in \figurename~\ref{plot2d}. We see from the figures that the same patterns are caught.

In case of computational efficiency, the ROM approximations are obtained by a speed-up factor $20$ over the FOM, \tablename~\ref{tab1}.
According to the energy criteria \eqref{cr2}, it requires global factor matrices of column size (number of modes) only $8-11$. The computed time averaged relative errors defined in \eqref{relerr} scales with $10^{-2}$. Detailed results for the number of modes and errors are presented in \tablename~\ref{tab2}.

\subsection{Brusselator equation}

We consider the three-dimensional Brusselator model \cite{Lin14}
in  the cubic domain $\Omega = [0,20]^3\subset\mathbb{R}^3$. We solve the problem through the matrix equation \eqref{sylvester31} with the mesh sizes $h_x=h_y=h_z=0.667$ for the number of grid points $n_1=n_2=n_3=31$. The final time is taken as $t_f=10$ with the time step size $\Delta t=0.01$, leading to the fourth dimension $n_4=1001$ of the snapshot tensors related to the time. For this problem, we fix the system parameters  $d_u=0.4$, $d_v=2$, $d_{uv}=0.02$, $\alpha =6$, $\beta =1$, and vary now the parameter $\theta :=d_{vu}$ in the set of admissible values $\textsc{P}=[19,23]$. As the finite training set of parameter $\theta$, we take the values (including the boundary values) uniformly distributed on $\textsc{P}$ with the increment $1$, i.e., $\textsc{P}_T=\{19,20,21,22,23\}$ with the number of sample parameter values $n_p=5$.

\begin{figure}[htb!]
	\centering
	\includegraphics[width=0.8\textwidth]{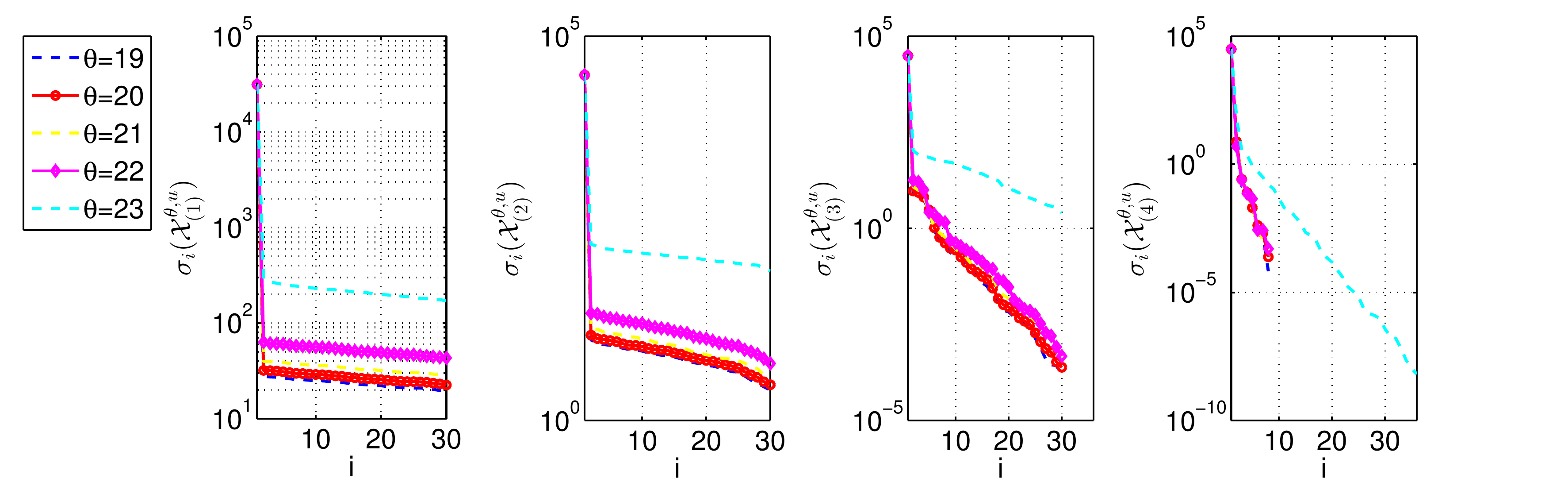}	
	\includegraphics[width=0.8\textwidth]{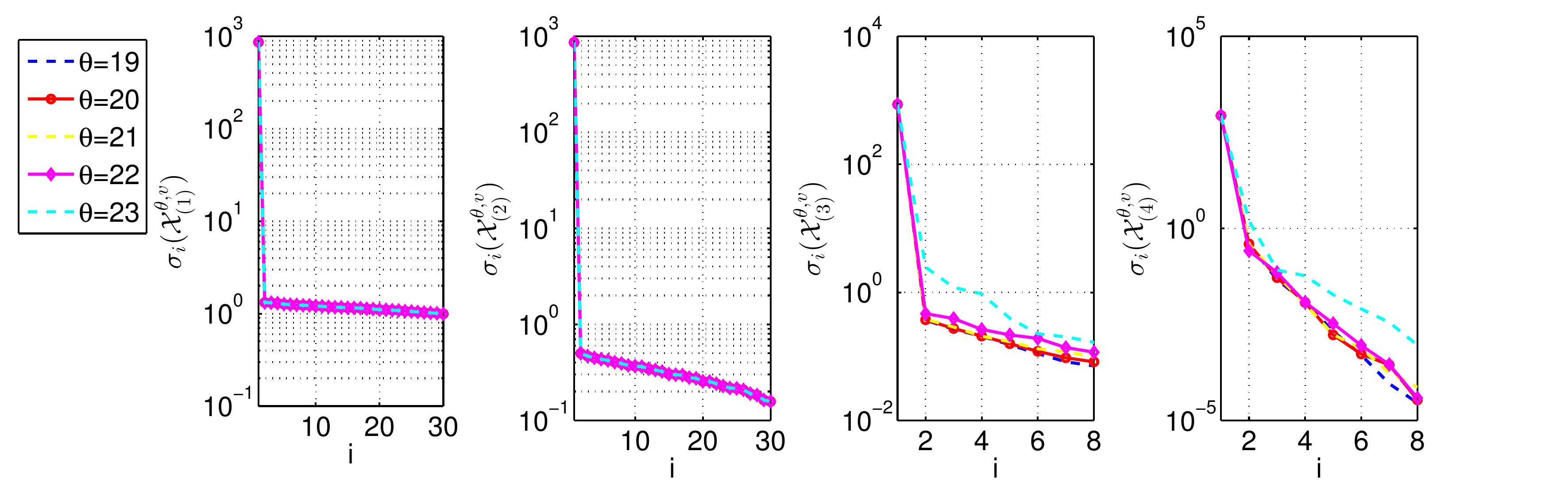}
	\caption{Brusselator model: Level I singular values of unfoldings ${\mathcal X}^{\theta,u}_{(j)}$ (top) and ${\mathcal X}^{\theta,v}_{(j)}$ (bottom) \label{sing3d}}
\end{figure}

\begin{figure}[htb!]
	\centering
	\includegraphics[width=0.6\textwidth]{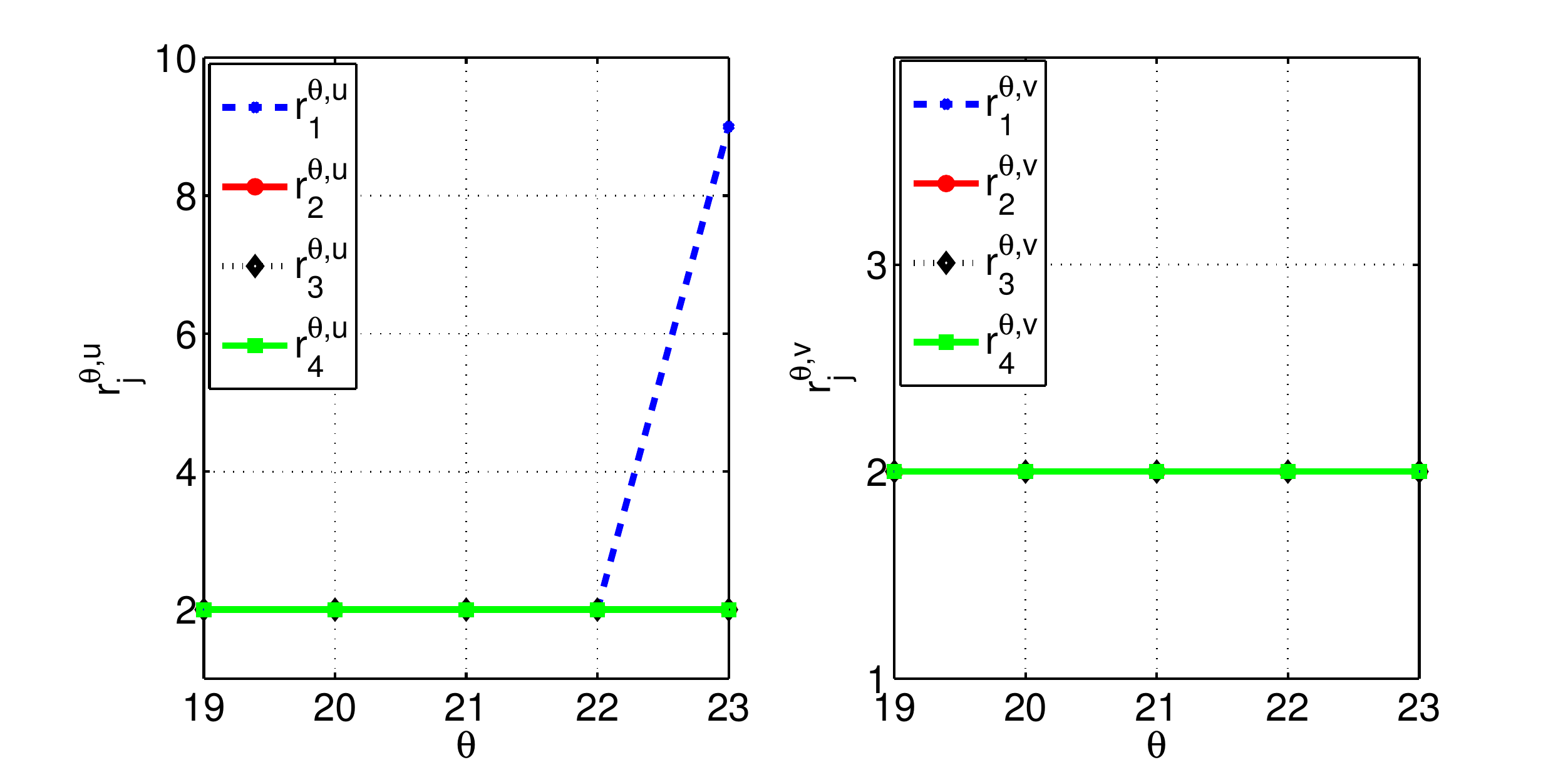}	
	\caption{Brusselator model: Sample parameter values vs target ranks for unfoldings ${\mathcal X}^{\theta,u}_{(j)}$ (left) and ${\mathcal X}^{\theta,v}_{(j)}$ (right) at Level I \label{numr3d}}
\end{figure}

In \figurename~\ref{sing3d}, we give the decay of the singular values $\sigma_i\left({\mathcal X}^{\theta,u}_{(j)}\right)$ and $\sigma_i\left({\mathcal X}^{\theta,v}_{(j)}\right)$ of the unfoldings ${\mathcal X}^{\theta,u}_{(j)}$ and ${\mathcal X}^{\theta,v}_{(j)}$ of the order-$4$ snapshot tensors ${\mathcal X}^{\theta,u}$ and ${\mathcal X}^{\theta,v}$, respectively, $j=1,2,3,4$, related to each sample parameter value $\theta\in\textsc{P}_T$. According to the criteria \eqref{cr1}, the computed target ranks $r^{\theta,u}_j$ and $r^{\theta,v}_j$ required by the HOSVD in the first level of nested POD are presented in \figurename~\ref{numr3d}. Similar to the previous example, it again shows that enough energetic part of the unfoldings are recovered.

\begin{figure}[htb!]
	\centering
	\includegraphics[width=0.7\textwidth]{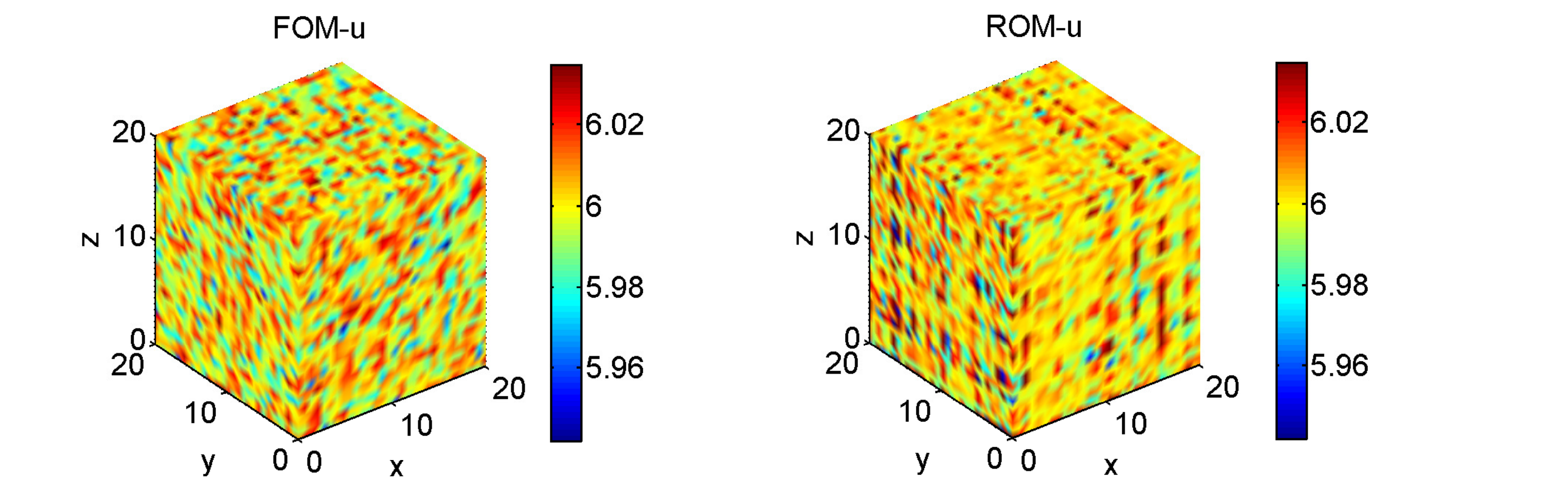}	
	\includegraphics[width=0.7\textwidth]{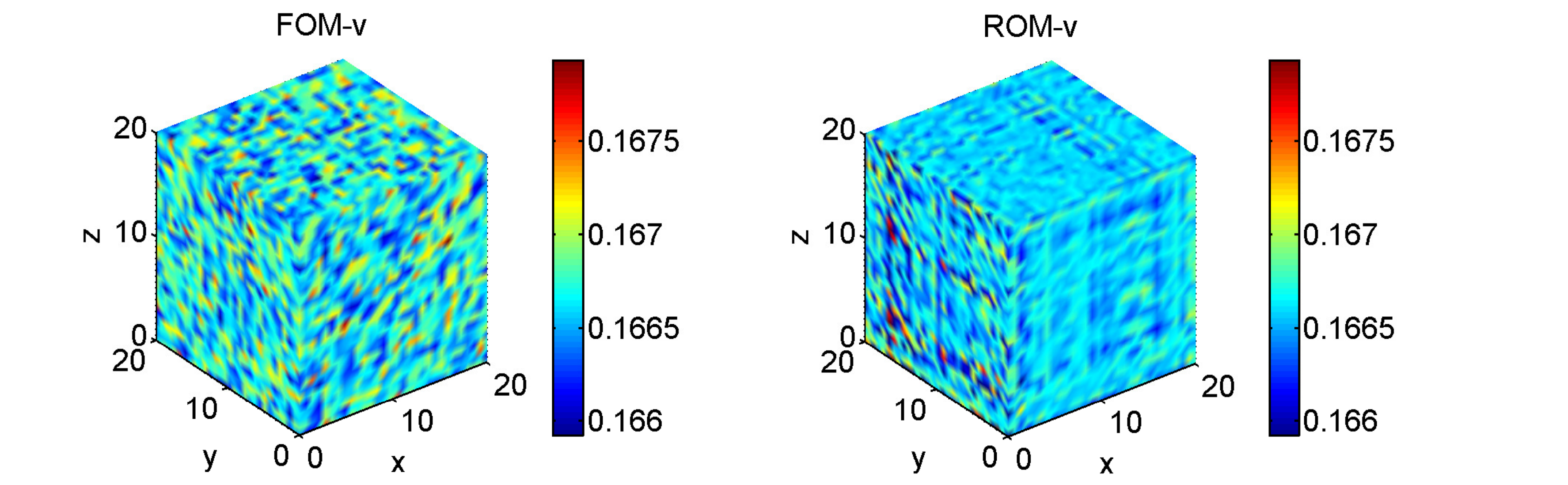}
	\caption{Brusselator model: FOM and ROM profiles at the final time for $\theta =21.5$ \label{plot3d}}
\end{figure}

The FOM solutions $u({\bm x}, t; \theta)$ and $v({\bm x}, t; \theta)$ together with the nonintrusive ROM approximations $\widehat{u}({\bm x}, t; \theta)$ and $\widehat{v}({\bm x}, t; \theta)$ at the final time $t_f=10$ for the parameter value $\theta =21.5\notin\textsc{P}_T$ are given in \figurename~\ref{plot3d}, where it can be seen that enough similar patterns are obtained.

In case of computational efficiency, \tablename~\ref{tab1} shows that the ROM approximations are obtained by a much greater speed-up factor, $130$, over the FOM compared with the speed-up factor obtained for the two-dimensional Schnakenberg model.
According to the energy criteria \eqref{cr2}, it requires global factor matrices of column size only $6-14$. The computed time averaged relative errors defined in \eqref{relerr} scales with $10^{-3}$. The detailed results for the Brusselator model can also be found in \tablename~\ref{tab2}.

\begin{table}[htb!]
\caption{Wall clock time (in seconds) and speed-up factors \label{tab1}}
\centering
\begin{tabular}{|lllrr|}
\hline%
 &  & & \textbf{Wall Clock Time} & \textbf{Speed-up}\\
\cline{4-5}%
\multirow{6}{*}{Schnakenberg}  & \multirow{4}{*}{Offline}  & FOMs  & 53.70  & \\
 & & Level I HOSVD Modes  & 223.29  & \\
	& & Level II POD Modes  & 0.03  & \\
	& & RBF Coefficients  & 0.93  &   \\
\cline{2-5}
& \multirow{2}{*}{Online $(\theta =0.65$)} & FOM  & 13.85  &   \\
	&	& ROM  & 0.70  &  {\bf 19.7}  \\
\hline
 \multirow{6}{*}{Brusselator} & \multirow{4}{*}{Offline}  & FOMs  & 586.77  &  \\
 & & Level I HOSVD Modes  & 74.70  &  \\
	& & Level II POD Modes  & 0.01  & \\
&	& RBF Coefficients  & 0.93  &  \\
\cline{2-5}
	&	\multirow{2}{*}{Online $(\theta =21.5)$} & FOM  & 123.30  &   \\
		&	 & ROM  & 0.95  &   {\bf 129.8}  \\
\hline
\end{tabular}
\end{table}

\begin{table}[htb!]
\caption{Time averaged relative errors and memory savings of compression \label{tab2}}
\centering
\begin{tabular}{|llrrcc|}
\hline%
 & &  \textbf{\#Modes ($u$, $v$)} &  \textbf{$\|u-\widehat{u}\|_{\text{rel}}$} & \textbf{$\|v-\widehat{v}\|_{\text{rel}}$} & \textbf{Saved Memory in \% }\\
\cline{3-6}%
\multirow{3}{*}{Schnakenberg $(\theta =0.65$)}  & $x$-direction ($\widehat{r}^{\cdot}_1$)  & 11, 10 & \multirow{3}{*}{9.54e-02} & \multirow{3}{*}{8.19e-02} & \multirow{3}{*}{\%99} \\
  & $y$-direction ($\widehat{r}^{\cdot}_2$) & 10, 9  &  & & \\
	& $t$-direction ($\widehat{r}^{\cdot}_3$) & 9, 8  &  & & \\
\hline
\multirow{4}{*}{Brusselator $(\theta =21.5)$}  & $x$-direction ($\widehat{r}^{\cdot}_1$)  & 14, 7 & \multirow{4}{*}{7.07e-03} & \multirow{4}{*}{6.63e-03} &  \multirow{4}{*}{\%99} \\
  & $y$-direction ($\widehat{r}^{\cdot}_2$) & 7, 7  &  & & \\
	& $z$-direction ($\widehat{r}^{\cdot}_3$) & 7, 7  &  & & \\
	& $t$-direction ($\widehat{r}^{\cdot}_4$) & 6, 6  &  & & \\
\hline
\end{tabular}
\end{table}

\begin{figure}[htb!]
	\centering
	\includegraphics[width=0.4\textwidth]{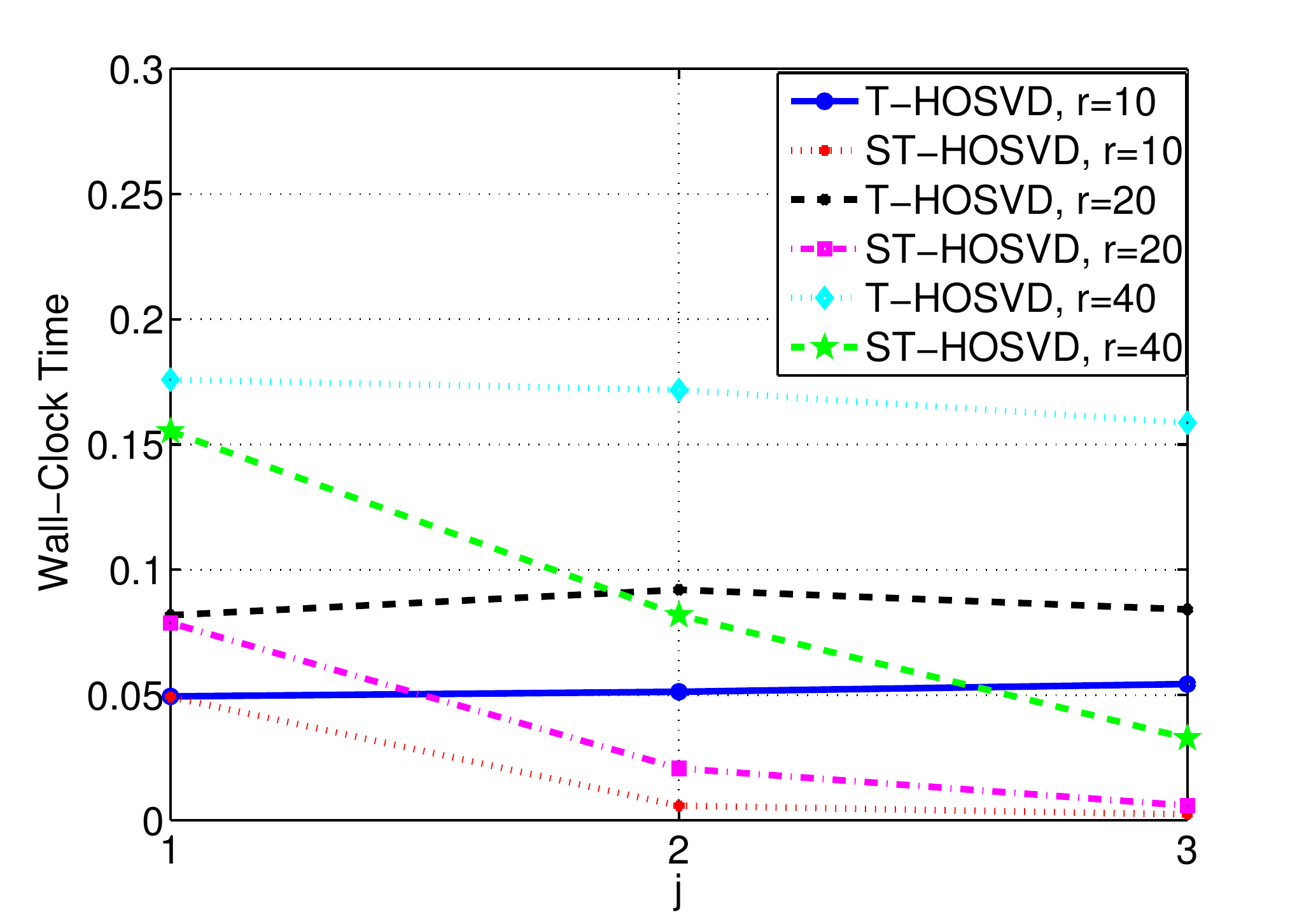}	
	\includegraphics[width=0.4\textwidth]{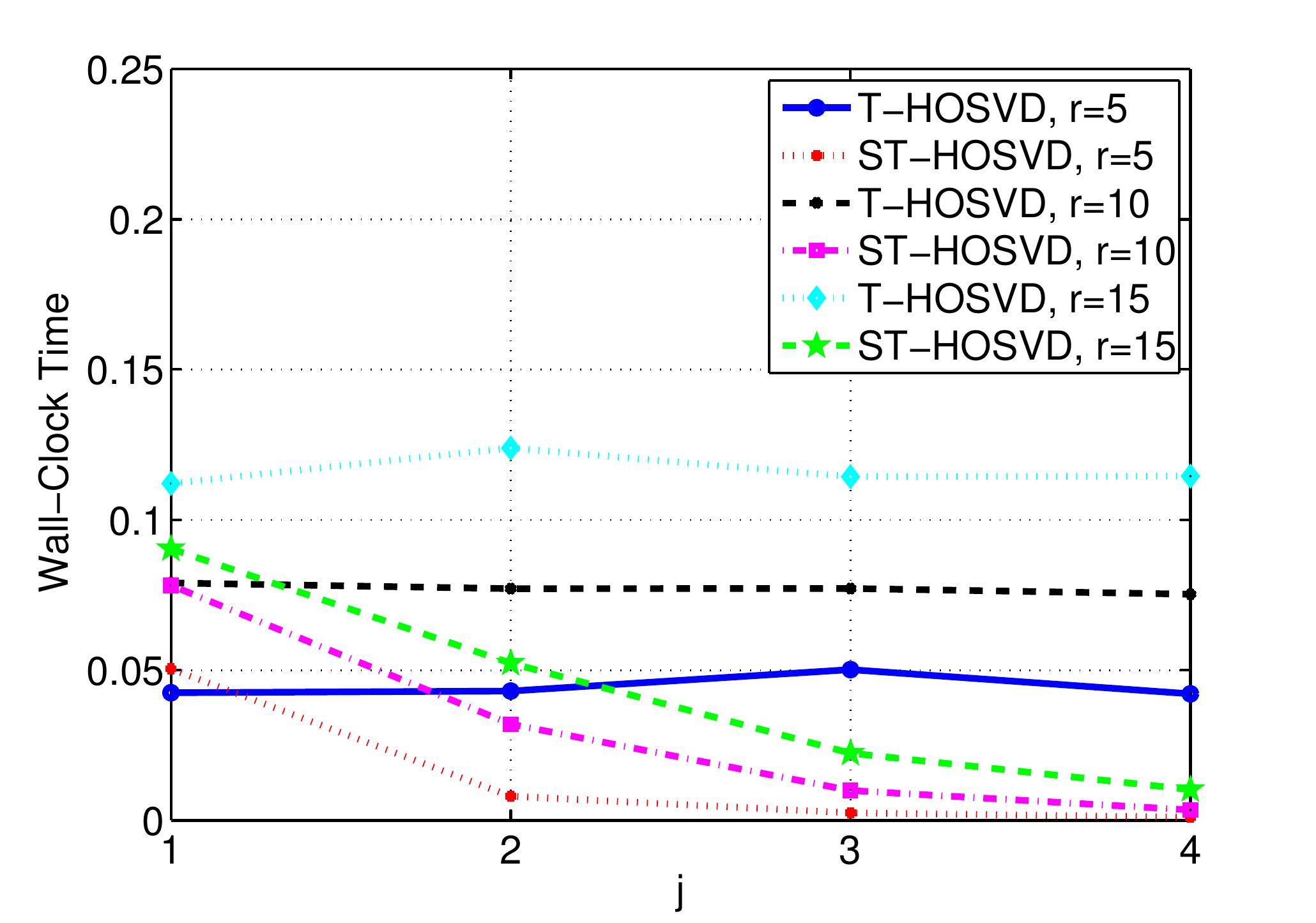}
	\caption{Wall-clock times of SVD computations elapsed for each unfolding ${\mathcal X}^{\theta_1,u}_{(j)}$ of the snapshot tensor ${\mathcal X}^{\theta_1,u}$ of Schnakenberg model (left) and Brusselator model (right), inside T-HOSVD and ST-HOSVD with different target rank $r$ \label{trunc_cpu}}
\end{figure}

We finally report the computational efficiency of the ST-HOSVD over T-HOSVD. To do this, we consider the order-$3$ and order-$4$ snapshot tensors ${\mathcal X}^{\theta_1,u}$ related to the component $u$ of both the two-dimensional Schnakenberg and three-dimensional Brusselator models with the same problem data considered above.
We apply T-HOSVD and ST-HOSVD to the snapshot tensors ${\mathcal X}^{\theta_1,u}$ with different values of target ranks $r:=r^{\theta_1,u}_1=\cdots =r^{\theta_1,u}_{d+1}$, and with the processing order $\bm{p}=[1,\ldots,d+1]$ for the ST-HOSVD. In \figurename~\ref{trunc_cpu}, we give the wall-clock times elapsed to make the SVD computations for each unfolding ${\mathcal X}^{\theta_1,u}_{(j)}$ of the order-$(d+1)$ snapshot tensors ${\mathcal X}^{\theta_1,u}$,  $d=2,3$.
It is clear that ST-HOSVD provides, in total, better computational efficiency compared to the T-HOSVD. Moreover, as the algorithm progresses, the time needed for the SVD computation decreases in ST-HOSVD, while it remains almost the same in the case of T-HOSVD.
This is because in the ST-HOSVD the dimension reduction is done sequentially, where the unfoldings are always obtained from the same full-rank snapshot tensor in T-HOSVD.

\section{Conclusions}
\label{sec:conc}

In this paper, we have developed nonintrusive ROMs exploiting the matrix/tensor based discretization of cross-diffusion systems in form of semilinear PDEs. The two-level approach for the construction of reduced bases through tensor decompositions with HOSVD instead of the classical SVD yields the reduced modes and reduced coefficients directly without necessitating further computation in the case of the vector based discretization. Numerical experiments with two-dimensional and three-dimensional cross-diffusion systems demonstrate the computational efficiency of the ROMs and the accuracy of the spatiotemporal patterns for new parameter values.



\end{document}